\documentclass{gtpart}     
\usepackage{pinlabel}  
\usepackage{graphicx}  
\title[Generalized Baumslag--Solitar groups]{On the isomorphism problem
for generalized \\Baumslag--Solitar groups}

\author{Matt Clay}
\givenname{Matt}
\surname{Clay}
\address{Mathematics Department\\
        University of Oklahoma\\\newline
        Norman, OK 73019\\ 
        USA}
\email{mclay@math.ou.edu}
\urladdr{}

\author{Max Forester}
\givenname{Max}
\surname{Forester}
\email{forester@math.ou.edu}
\urladdr{}

\arxivreference{0710.2108}
\arxivpassword{dzapc}

\volumenumber{}
\issuenumber{}
\publicationyear{}
\papernumber{}
\startpage{}
\endpage{}
\doi{}
\MR{}
\Zbl{}
\received{}
\revised{}
\accepted{}
\published{}
\publishedonline{}
\proposed{}
\seconded{}
\corresponding{}
\editor{}
\version{}

\theoremstyle{plain}
\newtheorem{theorem}{Theorem}
\newtheorem{proposition}{Proposition}[section]
\newtheorem{corollary}{Corollary}[section]
\newtheorem{lemma}{Lemma}[section]
\makeatletter
  \let\c@proposition\c@theorem
  \let\c@corollary\c@theorem
  \let\c@lemma\c@theorem
\makeatother

\theoremstyle{definition}
\newtheorem{definition}{Definition}[section]
\newtheorem{remark}{Remark}[section]
\newtheorem{example}{Example}[section]
\newtheorem{notation}{Notation}[section]
\makeatletter
  \let\c@definition\c@theorem
  \let\c@remark\c@theorem
  \let\c@example\c@theorem
  \let\c@notation\c@theorem
\makeatother

\numberwithin{theorem}{section}
\renewcommand{\theenumi}{\alph{enumi}}

\newcommand{\zindex}[3]{\put(#1,#2){\makebox(0,0){${#3}$}}}
\newenvironment{pict}[2]%
	{\setlength{\unitlength}{1mm}
	\begin{center}
	\begin{picture}(#1,#2)
	\scriptsize
}%
	{\end{picture}
	\end{center}
 	\noindent}

\DeclareFontFamily{OML}{rsfs}{\skewchar\font'177}
\DeclareFontShape{OML}{rsfs}{m}{n}{ <5> <6> rsfs5 <7> <8> <9> rsfs7
  <10> <10.95> <12> <14.4> <17.28> <20.74> <24.88> rsfs10 }{}
\DeclareMathAlphabet{\mathfs}{OML}{rsfs}{m}{n}

\newcommand{\As}{{\mathfs A}}

\newcommand{\DS}{\ensuremath{{\mathfs D}}} 
\newcommand{\LG}{\ensuremath{\mathsf{RLG}}} 
\newcommand{\pLG}{\ensuremath{\widehat{\LG}}} 
\newcommand{\Sl}{\ensuremath{\mathsf{S}}} 
\newcommand{\PLGI}{\ensuremath{\mathsf{P}}} 

\newcommand{\abs}[1]{\left| {#1} \right|}

\renewcommand{\geq}{\geqslant}
\renewcommand{\leq}{\leqslant}
\renewcommand{\epsilon}{\varepsilon}

\newcommand{\iv}{o} 
\newcommand{\tv}{t}

\DeclareMathOperator{\Out}{Out}

\begin{document}

\begin{abstract}    
Generalized Baumslag--Solitar groups (GBS groups) are groups that act
on trees with infinite cyclic edge and vertex stabilizers. Such an
action is described by a \emph{labeled graph} (essentially, the
quotient graph of groups). This paper addresses the problem of
determining whether two given labeled graphs define isomorphic groups;
this is the \emph{isomorphism problem} for GBS groups. There are two
main results and some applications. First, we find necessary and
sufficient conditions for a GBS group to be represented by only
finitely many reduced labeled graphs. These conditions can be checked
effectively from any labeled graph. Then we show that the isomorphism
problem is solvable for GBS groups whose labeled graphs have first Betti
number at most one. 
\end{abstract}

\maketitle


\section{Introduction}\label{sc:intro}

A \emph{generalized Baumslag--Solitar group} (or GBS group) is a group that
acts on a tree with infinite cyclic edge and vertex stabilizers.\footnote{In
  this paper we will only consider finitely generated GBS groups, so
  finite generation will be added to the definition; see
\fullref{ssc:gbs}.} The tree (together with the group action) is called a 
\emph{GBS tree}. A GBS tree can be described by a \emph{labeled graph},
which is a connected graph whose oriented edges are each labeled by a
non-zero integer. This information is enough to specify a graph of groups
encoding the GBS tree. 

A GBS group $G$ may have many labeled graph descriptions. Even if one
restricts to \emph{reduced} labeled graphs, which are in some sense the
simplest ones, there may be infinitely many distinct such graphs defining
$G$. It can also happen that there is only one reduced graph, or finitely
many. In these latter cases, useful information about $\Out(G)$ can be
obtained, as in Gilbert--Howie--Metaftsis--Raptis \cite{ar:GHMR00}, Pettet
\cite{ar:P99}, and Levitt \cite{ar:L07}. Other aspects
of GBS groups have been studied by Kropholler, Whyte, Levitt, and
others. See Kropholler \cite{ar:K90,ar:K90b}, Whyte \cite{ar:W01}, Clay
\cite{ar:Cpp}, and Forester \cite{ar:F03,ar:F06}
for details on various algebraic and geometric properties of GBS groups.

The variety of labeled graph descriptions of GBS groups is partly what
makes them interesting. For instance, they demonstrate the extent to
which JSJ decompositions of groups can fail to be unique. On the other
hand, this variety can also be a source of difficulty, such as when
studying automorphisms. A given labeled graph need not be invariant, for
instance. Even the basic problem of recognizing a given GBS group from
one of its labeled graphs is not at all clear. 

The \emph{isomorphism problem} for GBS groups is the problem of
determining algorithmically whether two given labeled graphs define
isomorphic GBS groups. This problem has only been shown to be solvable in
limited special cases. It is trivially solvable for the \emph{rigid} GBS
groups, which are those having a unique reduced labeled graph. These 
groups were characterized by Levitt \cite{ar:L05a} (see also 
Gilbert et al. \cite{ar:GHMR00}, Pettet \cite{ar:P99},
and Forester \cite{ar:F02}). 

Levitt showed that the isomorphism problem is solvable in the case of GBS
groups $G$ such that $\Out(G)$ does not contain a non-abelian free group
\cite{ar:L07}. He also solved the isomorphism problem for
$2$--generator GBS groups \cite{ar:Lprep}. Both of these results rely on
having an explicit characterization of the class of groups being
considered. 

In \cite{ar:F06} the isomorphism problem was solved for GBS
groups whose 
modular groups contain no integers other than $\pm 1$. Equivalently (see
Levitt \cite{ar:L07}) these are the GBS groups not containing any solvable
Baumslag--Solitar group $BS(1,n)$ with $n > 1$. It is worth recalling the
main steps of the proof. First it was shown that any two such graphs are
related by slide moves, without leaving the set of reduced graphs. Then
it was shown that such a group is represented by only finitely many
reduced labeled graphs. Thus, this set can be searched and enumerated
effectively, and membership is decidable. 

For the general isomorphism problem, it is useful to understand the space
of reduced labeled graphs related to a given one by sequences of slide
moves. We want to know whether this space is infinite, and whether it
includes all reduced labeled graphs for the given group. 
To this end, there is a property of edges that plays a key role: edges can 
be \emph{mobile} or \emph{non-mobile} (see \fullref{def:mobile}). One of
our main technical results  
is \fullref{co:finitefirst}, which shows that in any sequence of
slide moves, the non-mobile edges may be slid first and one at a
time. From this we deduce information on the slide space of a labeled
graph, including our first main result: 

\begin{theorem}\label{th:finite}
Let $G$ be a GBS group other than $BS(1,n)$, represented by a reduced
labeled graph $\Gamma$. Then $G$ has finite reduced labeled graph space 
if and only if $\Gamma$ has no mobile edges. 
\end{theorem}

We also show that mobility of edges can be tested algorithmically
(\fullref{rm:mobile}), so the property of the theorem is decidable. (The
case of $BS(1,n)$ is clear as well: the reduced labeled graph is unique if $n
\not= -1$, and $BS(1,-1)$ has two reduced labeled graphs.) 
One consequence of \fullref{th:finite} is \fullref{th:iso-finite}, which
solves the isomorphism problem in the case where one labeled graph has no
mobile edges. 

Next we consider the case of GBS groups whose labeled graphs have first
Betti number one. (The Betti number zero case is covered by
Forester \cite{ar:F06}.) The primary goal of the rest of the paper is the
following theorem:  

\begin{theorem} \label{th:iso} 
There is an algorithm which, given two labeled graphs, one of which has
first Betti number at most one, determines whether the two GBS groups are
isomorphic. 
\end{theorem}

These are two cases, which behave rather differently: the ascending case
and the non-ascending case. In the ascending case, there is a structure
theorem (\fullref{th:asc-structure}) which says that the group is
uniquely determined by certain invariants, which can be computed by
putting the labeled graph into a normal form. These invariants are
defined and proved invariant with the aid of Theorem 1.1 of
Clay--Forester \cite{ar:CF}, which shows that any two reduced labeled
graphs are related 
by slide, induction, and $\As^{\pm 1}$--moves between reduced labeled
graphs (see \fullref{ssc:gbs} for the definitions of these moves). 

The non-ascending case is somewhat simpler, since any two reduced
labeled graphs representing the same group are related by slide
moves. (In particular, one may keep track of individual edges.) However,
even though we can define normal forms, they are much less rigid than in
the ascending case. For instance, there is no canonical edge with which
to compare other edges, unlike ascending normal forms. 

We show that given $G$, there are only finitely many reduced labeled
graphs in normal form, and these can be enumerated effectively. The
solution to the isomorphism problem is then similar to the case proved in
\cite{ar:F06}.

\subsection*{Acknowledgements} The second author is partially supported
by NSF grant \ DMS-0605137.

\section{Preliminaries}\label{sc:prelim}

\subsection{Deformation spaces}\label{ssc:def}
A graph $\Gamma$ is given by $(V(\Gamma),E(\Gamma),\iv,\tv,\bar{
\ })$ where $V(\Gamma)$ are the vertices, $E(\Gamma)$ are the oriented
edges, $\iv,\tv \co E(\Gamma) \to V(\Gamma)$ are the originating and
terminal vertex maps and $\bar{ \ }\co E(\Gamma) \to E(\Gamma)$ is a
fixed point free involution, which reverses the orientations of edges.
An \emph{edge path} $\gamma = (e_{0},\ldots,e_{k})$ is a sequence of
edges such that $\tv(e_{i}) = \iv(e_{i+1})$ for $i = 0,\ldots,k-1$.  A
\emph{loop} is an edge $e \in E(\Gamma)$ such that $\iv(e) = \tv(e)$. A
\emph{geometric   edge} is a pair of the form $\{e, \bar{e}\}$. 

Let $G$ be a group.  A \emph{$G$--tree} is a simplicial tree $T$ together
with an action of $G$ by simplicial automorphisms, without inversions 
(that is, $ge \neq \bar{e}$ for all $g \in G, e \in E(T)$).  Two
$G$--trees are considered equivalent if there is a $G$--equivariant
isomorphism between them. The quotient graph $T/G$ has the structure of a
graph of groups with a marking (an identification of $G$ with the
fundamental group of the graph of groups). 

Given a $G$--tree $T$, a subgroup $H \subseteq G$ is \emph{elliptic}
if it fixes a point of $T$.  There are two moves one can perform
on a $G$--tree without changing the elliptic subgroups, called
\emph{collapse and expansion moves}; they correspond to the natural
isomorphism $A \ast_{B} B \cong A$. The exact definition is as follows.

\begin{definition}\label{def:elementarymoves}
An edge $e$ in a $G$--tree $T$ is \emph{collapsible} if $G_{e} =
G_{\iv(e)}$ and its endpoints are not in the same orbit. If one collapses
$\{e, \bar{e}\}$ and all of its translates to vertices, the resulting
$G$--tree is said to be obtained from $T$ by a \emph{collapse move}. The
reverse of this move is called an \emph{expansion move}. 

A $G$--tree is \emph{reduced} if it does not admit a
collapse move. An \emph{elementary deformation} is a finite sequence of
collapse and expansion moves. Given a $G$--tree $T$, the \emph{deformation
  space} $\DS$ of $T$ is the set of all $G$--trees related to $T$ by an
elementary deformation. If $T$ is cocompact then $\DS$ is equivalently the
set of all $G$--trees having the same elliptic subgroups as $T$, by
Forester \cite{ar:F02}. 
\end{definition}

\subsection{Generalized Baumslag--Solitar groups}\label{ssc:gbs}

A group $G$ that acts on a tree with infinite cyclic
stabilizers is called a \emph{generalized Baumslag--Solitar group} (or
\emph{GBS group}). In this paper, for simplicity, we also require $G$ to
be finitely generated (this convention is not followed in
\cite{ar:F02,ar:F03,ar:F06}). 
The tree is called a \emph{GBS tree}. In the quotient
graph of groups, every vertex and edge group is isomorphic to $\Z$,
and each inclusion map $G_{e} \hookrightarrow G_{\iv(e)}$ is given
by multiplication by a non-zero integer.  This data can be effectively
represented in a \emph{labeled graph}.  Specifically, a labeled graph
is a pair $(\Gamma,\lambda)$ where $\Gamma$ is a finite connected graph
and $\lambda\co E(\Gamma) \to \Z - \{0\}$ is a function, called
the labeling.  Given a choice of generators of $G_e$ and $G_{\iv(e)}$,
the map $G_{e} \hookrightarrow G_{\iv(e)}$ is multiplication by
$\lambda(e)$. Replacing a generator of an edge group $G_{e}$ by its
inverse interchanges the signs of $\lambda(e)$ and $\lambda(\bar{e})$;
replacing a generator of a vertex group $G_{v}$ by its inverse
interchanges the signs of $\lambda(e)$ for all edges $e$ originating
at $v$.  These operations are called \emph{admissible sign changes}. 
This is the only ambiguity in the labels of a labeled graph. We will
sometimes refer to $(\Gamma, \lambda)$ simply as $\Gamma$. 

A $G$--tree is
\emph{elementary} if there is a $G$--invariant point or line, and is
\emph{non-elementary} otherwise. By Forester \cite[Lemma~2.6]{ar:F03}, a
GBS tree 
is elementary if and only if the group is isomorphic to $\Z$, $\Z\times
\Z$, or the Klein bottle group. Thus a GBS group not isomorphic to one of
these three groups is called a \emph{non-elementary GBS group}. 

In a non-elementary GBS group, the elliptic subgroups arising from any
GBS tree are characterized algebraically \cite[Lemma~2.5]{ar:F03}.
Therefore, any two such $G$--trees lie in the same deformation space.
In particular, any two labeled graphs representing the same
non-elementary group are related by an elementary deformation. 
Whenever we speak of a deformation space for a non-elementary GBS
group, we will always be referring to this canonical deformation space.
For a description of this canonical deformation space associated to
the classical Baumslag--Solitar groups $BS(p,q)$, see Clay \cite{ar:Cpp}.
Unless otherwise stated, all GBS groups considered here will be assumed
to be non-elementary.

In a labeled graph, a loop $e$
with label $\pm 1$ is called an \emph{ascending loop}. It is a
\emph{strict ascending loop} if $\lambda(\bar{e}) \not= \pm 1$. 
A loop $e$ is a \emph{virtual ascending loop} if $\lambda(e)$ divides
$\lambda(\bar{e})$, and is a \emph{strict virtual ascending loop} if, in
addition, $\lambda(\bar{e}) \not= \pm \lambda(e)$. 
A GBS deformation space is \emph{ascending} if it contains a 
GBS tree whose labeled graph has a strict ascending loop. We also say
that $G$ is ascending. Otherwise the deformation space (or the group) is
called \emph{non-ascending}. 

Now we define various moves between GBS trees, all of which are
elementary deformations. The moves in \fullref{def:oldmoves} are
discussed more fully (in the general setting of $G$--trees) in
\cite{ar:F06}. In particular, slides and inductions can be factored as an
expansion followed by a collapse. A general discussion of $\As^{\pm 
1}$--moves can be found in Clay--Forester \cite{ar:CF}. 

In the diagrams below, each label $\lambda(e)$ is pictured next to the
endpoint $\iv(e)$. We begin with the elementary moves, which look
as follows (modulo admissible sign changes):

\begin{pict}{90}{10}
\thicklines
\put(10,5){\circle*{1}}
\put(25,5){\circle*{1}}

\put(10,5){\line(1,0){15}}

\thinlines
\put(45,6.5){\vector(1,0){15}}
\put(60,2.2){\vector(-1,0){15}}
\zindex{52.5}{8.2}{\mbox{collapse}}
\zindex{52.5}{4}{\mbox{expansion}}

\put(25,5){\line(3,5){3}}
\put(25,5){\line(3,-5){3}}

\put(10,5){\line(-5,3){5}}
\put(10,5){\line(-5,-3){5}}

\put(80,5){\circle*{1}}
\put(80,5){\line(-5,3){5}}
\put(80,5){\line(-5,-3){5}}

\put(80,5){\line(3,5){3}}
\put(80,5){\line(3,-5){3}}

\scriptsize
\zindex{8.5}{8}{a}
\zindex{8.5}{2}{b}
\zindex{12}{6.5}{n}
\zindex{23}{6.5}{1}
\zindex{28.5}{7.5}{c}
\zindex{28.5}{2.5}{d}

\zindex{78.5}{8}{a}
\zindex{78.5}{2}{b}
\zindex{84.5}{7.5}{nc}
\zindex{84.5}{2.5}{nd}
\end{pict}
Thus a GBS tree is reduced if and only if its labeled graph does not
contain an edge with distinct endpoints and label $\pm 1$. 

\begin{definition}\label{def:oldmoves}
A \emph{slide move} between GBS trees takes one of the following two
forms: 

\begin{pict}{100}{11}
\thicklines
\put(75,3){\circle*{1}}
\put(90,3){\circle*{1}}
\put(75,3){\line(1,0){15}}
\put(90,3){\line(-1,2){4}}

\thinlines
\put(47.5,3){\vector(1,0){10}}
\zindex{52.5}{5}{\mbox{slide}}
\put(90,3){\line(5,3){5}}
\put(90,3){\line(5,-3){5}}
\put(75,3){\line(-5,3){5}}
\put(69,3){\line(1,0){6}}
\put(75,3){\line(-5,-3){5}}

\scriptsize
\zindex{77}{1.5}{m}
\zindex{88}{1.5}{n}
\zindex{86.5}{6}{\ell n}

\thicklines
\put(10,3){\circle*{1}}
\put(25,3){\circle*{1}}
\put(10,3){\line(1,0){15}}
\put(10,3){\line(1,2){4}}

\thinlines
\put(25,3){\line(5,3){5}}
\put(25,3){\line(5,-3){5}}
\put(10,3){\line(-5,3){5}}
\put(4,3){\line(1,0){6}}
\put(10,3){\line(-5,-3){5}}

\scriptsize
\zindex{12}{1.5}{m}
\zindex{23}{1.5}{n}
\zindex{9.5}{7}{\ell m}
\end{pict} 
or

\begin{pict}{100}{10}
\thicklines
\put(82.5,5){\circle*{1}}
\put(87.5,5){\circle{10}}
\put(72.5,5){\line(1,0){10}}

\thinlines
\put(47.5,5){\vector(1,0){10}}
\zindex{52.5}{7}{\mbox{slide}}
\put(82.5,5){\line(-5,-3){4.5}}
\put(82.5,5){\line(-1,-4){1.2}}
\put(82.5,5){\line(-1,4){1.2}}

\scriptsize
\zindex{85.1}{3.5}{m}
\zindex{84.7}{6.5}{n}
\zindex{79.9}{6.8}{\ell n}

\thicklines
\put(17.5,5){\circle*{1}}
\put(22.5,5){\circle{10}}
\put(7.5,5){\line(1,0){10}}

\thinlines
\put(17.5,5){\line(-5,-3){4.5}}
\put(17.5,5){\line(-1,-4){1.2}}
\put(17.5,5){\line(-1,4){1.2}}

\scriptsize
\zindex{20.1}{3.5}{m}
\zindex{19.7}{6.5}{n}
\zindex{14.5}{6.8}{\ell m}
\end{pict} 
An \emph{induction move} between GBS trees is as follows:

\smallskip
\begin{pict}{80}{10}
\thicklines
\put(6,5){\circle{10}}
\put(11,5){\circle*{1}}
\put(68,5){\circle{10}}
\put(73,5){\circle*{1}}

\thinlines
\put(11,5){\line(1,1){4}}
\put(11,5){\line(1,-1){4}}
\put(73,5){\line(1,1){4}}
\put(73,5){\line(1,-1){4}}

\scriptsize
\zindex{15}{7}{a}
\zindex{15}{3.5}{b}
\zindex{9}{6.7}{1}
\zindex{8}{3.3}{\ell m}
\put(31.5,5){\vector(1,0){17}}
\put(48.5,5){\vector(-1,0){17}}
\zindex{40}{7}{\mbox{induction}}
\zindex{78}{7}{\ell a}
\zindex{78}{3.3}{\ell b}
\zindex{71}{6.7}{1}
\zindex{70}{3.3}{\ell m}
\end{pict}
Both directions of the move are considered induction moves. This move
decomposes into an elementary deformation as follows: 

\begin{pict}{100}{10}
\thicklines
\put(48,5){\oval(10,10)[b]}

\put(88,5){\circle{10}}
\put(93,5){\circle*{1}}

\thinlines
\put(6,5){\circle{10}}
\put(11,5){\circle*{1}}

\put(48,5){\circle{10}}
\put(43,5){\circle*{1}}
\put(53,5){\circle*{1}}

\put(11,5){\line(1,1){4}}
\put(11,5){\line(1,-1){4}}

\put(53,5){\line(1,1){4}}
\put(53,5){\line(1,-1){4}}

\put(93,5){\line(1,1){4}}
\put(93,5){\line(1,-1){4}}

\scriptsize
\zindex{15}{7}{a}
\zindex{15}{3.2}{b}
\zindex{9}{6.7}{1}
\zindex{8}{3.3}{\ell m}

\put(22,5){\vector(1,0){12}}
\zindex{28}{6.5}{\mbox{exp. }}

\zindex{57}{7}{a}
\zindex{57}{3.2}{b}
\zindex{51}{6.7}{1}
\zindex{50.5}{3}{m}
\zindex{41.5}{6.7}{\ell}
\zindex{41.5}{3}{1}

\put(64,5){\vector(1,0){12}}
\zindex{70}{7}{\mbox{coll. }}

\zindex{98}{7}{\ell a}
\zindex{98}{3}{\ell b}
\zindex{91}{6.7}{1}
\zindex{90}{3.3}{\ell m}
\end{pict}%
\end{definition}

\begin{definition}\label{def:A}
Next we discuss \emph{$\As^{\pm 1}$--moves}, defined in \cite{ar:CF}. An 
$\As^{-1}$--move is an induction followed by a collapse, with the
following description. It is required that $k, \ell \not= \pm 1$, and
that the left hand vertex has no other edges incident to it. 

\smallskip
\begin{pict}{90}{10}
\thicklines
\put(6,5){\circle{10}}
\put(11,5){\circle*{1}}
\put(73,5){\circle{10}}
\put(78,5){\circle*{1}}
\put(24,5){\circle*{1}}
\put(11,5){\line(1,0){13}}

\thinlines
\put(24,5){\line(1,1){4}}
\put(24,5){\line(1,-1){4}}
\put(78,5){\line(1,1){4}}
\put(78,5){\line(1,-1){4}}

\scriptsize
\zindex{28}{7}{a}
\zindex{28}{3.5}{b}
\zindex{9}{6.7}{1}
\zindex{8}{3.3}{\ell m}
\zindex{13}{7}{\ell}
\zindex{22}{7}{k}
\put(40,6.5){\vector(1,0){15}}
\put(55,2.2){\vector(-1,0){15}}
\zindex{47.5}{8.5}{\As^{-1}}
\zindex{47.5}{4}{\As}

\zindex{82}{7}{a}
\zindex{82}{3.3}{b}
\zindex{75.5}{6.7}{k}
\zindex{74}{3.3}{k\ell m}
\end{pict}
The induction move changes the label $\ell$ to $1$, after which the edge
is collapsed. 

Note that before the move, the loop is a strict ascending loop,
and after, the loop is not ascending. Thus an $\As^{-1}$--move
removes an ascending loop, and its reverse, called an \emph{$\As$--move},
adds one. 
\end{definition}

\begin{remark}
$\As^{\pm 1}$--moves preserve the property of being reduced. The same is
not always true of slide or induction moves, unless one is in a
non-ascending deformation space. Also, an induction or $\As^{\pm 1}$--move
can only occur in an ascending deformation space. 
\end{remark}

We will make extensive use of the following result, which is the main
theorem of \cite{ar:CF}, and its corollary. 

\begin{theorem}\label{th:moves}
In a deformation space of cocompact $G$--trees, any two reduced trees are
related by a finite sequence of slides, inductions, and $\As^{\pm
  1}$--moves, with all intermediate trees reduced. 
\end{theorem}

\begin{corollary}\label{co:slides}
In a non-ascending deformation space of cocompact $G$--trees, any two
reduced trees 
are related by a finite sequence of slide moves, with all intermediate
trees reduced. Moreover, if $e$ is an edge of $T$ and a deformation from
$T$ to $T'$ never collapses $e$, then there is a sequence of slide moves
from $T$ to $T'$ in which no edge slides over $e$. 
\end{corollary}

The first statement of the corollary follows immediately from the
theorem, and has previously appeared as Forester \cite[Theorem
7.4]{ar:F06} and Guirardel--Levitt \cite[Theorem~7.2]{ar:GL07}. The
second statement is proved in \cite{ar:CF}.

\subsection{The modular homomorphism}\label{ssc:modular}

Let $G$ be a GBS group with labeled graph $(\Gamma, \lambda)$. 
There are two versions of the modular homomorphism $G \to
\Q^{\times}$, each with several descriptions; see Bass--Kulkarni 
\cite{ar:BK90}, Forester \cite{ar:F06}, and Kropholler \cite{ar:K90b}. In
this paper, it turns out to be more 
convenient to use the reciprocal of the usual definition, so we will
include this modification here. This makes it easier to keep track of
slide moves; see for example \fullref{def:e-edgepath} and 
\fullref{rm:e-edgepath}. We will mostly work with the \emph{signed modular 
  homomorphism} $q \co G \to \Q^{\times}$, defined as the composition $G
\to H_1(\Gamma) \to \Q^{\times}$ where the second map is given by 
\begin{equation} \label{eq:modular}
 (e_1, \ldots, e_k) \mapsto \prod_{i=1}^k
\frac{\lambda(\bar{e}_i)}{\lambda(e_i)}. 
\end{equation}
(The first map is given by first killing the normal closure of the
elliptic elements to obtain $\pi_1(\Gamma)$, and then abelianizing.) 
Equivalently, fix a non-trivial elliptic element $a \in G$. Then every
$g\in G$ satisfies a relation $g a^r g^{-1} = a^s$ in $G$ for
some non-zero integers $r$ and $s$, and the assignment $q(g) = s/r$ is a
well defined homomorphism, which agrees with the definition just given;
see Kropholler \cite{ar:K90b} or Levitt \cite{ar:L07}. 

The \emph{unsigned modular homomorphism} is simply $\abs{q}$, defined on
$H_1(\Gamma)$ by 
\[(e_1, \ldots, e_k) \mapsto \prod_{i=1}^k
\frac{\abs{\lambda(\bar{e}_i)}}{\abs{\lambda(e_i)}}. \] 
An equivalent definition is to choose any subgroup $V$ of $G$ 
commensurable with a vertex group, and assign to each $g \in G$ the 
positive rational number 
\[[V^g : V \cap V^g] \ / \ [V : V \cap V^g].\] 
See \cite{ar:F06} for a proof that this function agrees with $\abs{q(g)}$. 
We say that $(\Gamma,\lambda)$ is \emph{unimodular} if $\abs{q}$ is
trivial. 

Finally, there is also an \emph{orientation homomorphism} $G \to \{\pm
1\}$ defined by $g \mapsto q(g)/\abs{q(g)}$.  This homomorphism is also
defined on $H_1(\Gamma)$.  The next result shows that the GBS group
associated to a labeled graph is determined by the orientation
homomorphism and the absolute value of the labeling. Hence it often
suffices to consider \emph{positive labeled graphs}, i.e. labeled graphs
$(\Gamma,\lambda)$ such that $\lambda(e) > 0$ for all $e \in E(\Gamma)$. 

\begin{lemma}\label{lm:signs}
Let $\lambda$ and $\lambda'$ be labelings on a graph $\Gamma$ such that
$\abs{\lambda} = \abs{\lambda'}$. If their orientation homomorphisms
agree then $(\Gamma,\lambda)$ and $(\Gamma, \lambda')$ differ by
admissible sign changes. In particular, the corresponding GBS groups are
isomorphic. 
\end{lemma}

\begin{proof}
Let $\Omega \co H_1(\Gamma) \to \{\pm 1\}$ be the orientation homomorphism
of $(\Gamma,\lambda)$ and $(\Gamma,\lambda')$. Fix a maximal tree $T
\subseteq \Gamma$. Then every edge $e$ of $\Gamma - T$ determines a generator
$[e] \in H_1(\Gamma)$. 

By admissible sign changes, we can arrange that $\lambda$ and
$\lambda'$ agree, and are positive, on the edges of $T$. Then for any
edge $e$ in $\Gamma - T$ we have that $\Omega([e]) = 1$ if and only if
$\lambda(e)$ and $\lambda(\bar{e})$ have the same sign, if and only
if $\lambda'(e)$ and $\lambda'(\bar{e})$ have the same sign. Thus
$\lambda$ and $\lambda'$ can be made to agree on $e$ and $\bar{e}$
by an admissible sign change affecting $e, \bar{e}$ only. In this
way, $\lambda$ and $\lambda'$ can be made to agree on all of $\Gamma$. 
\end{proof}

\section{Labeled graph spaces}\label{sc:labeled}

From now on we consider only GBS groups and their canonical
deformation spaces.  Hence we will always refer to $G$ instead of this
deformation space.  

\begin{definition}\label{def:labeled}
For a GBS group $G$, we denote by $\LG(G)$ the set of 
reduced 
labeled graphs representing $G$. Let $\LG^+(G)$ be the set of positive
reduced labeled graphs representing $G$. Note that this latter set is
non-empty only if the orientation homomorphism is trivial. 
\end{definition}

Our goal in this section is to establish a criterion, which can 
be checked in terms of any labeled graph in $\LG(G)$, that characterizes
when $\LG(G)$ is finite.  Notice that if $G$ is ascending and $G \neq
BS(1,n)$, then $|\LG(G)| = \infty$.  If $G = BS(1,n)$ or $G = \Z$,
then $\LG(G)$ consists of a single point (unless $G = BS(1,-1)$, in
which case $|\LG(G)| = 2$).  Therefore, we are mainly concerned 
with determining when a non-ascending GBS group satisfies $|\LG(G)| =
\infty$. However, we will need to prove a more general statement, as we
do not have an algorithm to determine whether a given GBS group is
ascending.

\subsection{Monotone cycles and mobile edges}\label{ssc:mobile} 
If $(e_0, \ldots, e_n)$ is an edge path in $\Gamma$, we define $q(e_0,
\ldots, e_n)$ by formula \eqref{eq:modular}. This is also denoted
$q_{\Gamma}(e_0, \ldots, e_n)$. 
\begin{definition}\label{def:e-edgepath}
Let $\Gamma$ be a labeled graph for $G$ and $e \in E(\Gamma)$.  An 
edge path $(e_{0},\ldots,e_{n})$ is an $e$--\emph{edge path} if:

\begin{enumerate}
    
    \item $e_{i} \neq e, \bar{e}$ for $i = 0,\ldots,n$;
    
    \item $\iv(e) = \iv(e_{0})$; and
    
    \item $\lambda(e_{i})$ divides
    $\lambda(e)q(e_{0},\ldots,e_{i-1})$ for $i = 0,\ldots,n$.
    
\end{enumerate}
An $e$--edge path is an $e$--\emph{integer cycle} if, in addition we 
have:

\begin{enumerate}
\setcounter{enumi}{3}    
    \item $\iv(e_{0}) = \tv(e_{n})$; and
    
    \item $q(e_{0},\ldots,e_{n}) \in \Z$.
    
\end{enumerate}
If $|q(e_{0},\ldots,e_{n})| \neq 1$ we say that the
$e$--edge path or the $e$--integer cycle is \emph{strict}.
\end{definition}

\begin{remark}\label{rm:e-edgepath}
The first three conditions are necessary and sufficient to be able to
slide $e$ along $(e_0, \ldots, e_n)$. The resulting label on the edge $e$ is
$\lambda(e)q(e_{0},\ldots,e_{n})$.  Hence $e$ may slide
repeatedly along an $e$--integer cycle.  Also notice that any path
obtained by tightening an $e$--edge path (respectively, $e$--integer
cycle) is an $e$--edge path (respectively, $e$--integer cycle).
\end{remark}

\begin{definition}\label{def:monotone}
An edge path $(e_{0},\ldots,e_{n},e)$ is a \emph{monotone cycle} if
$(e_{0},\ldots,e_{n})$ is an $\bar{e}$--edge path and
$q(e_{0},\ldots,e_{n},e) \in \Z$.  An edge $e$ is a monotone cycle
if $e$ is a loop and $q(e) \in \Z$.  A monotone cycle is \emph{strict}
if the modulus is not equal to $\pm 1$.
\end{definition}

\begin{remark}\label{rm:monotone}
Suppose $(e_{0},\ldots,e_{n},e)$ is a monotone cycle.  Since
$(e_{0},\ldots,e_{n})$ is an $\bar{e}$--edge path, $\iv(e_{0}) =
\iv(\bar{e}) = \tv(e)$.  Hence a monotone cycle is a cycle.
Further, since $\bar{e} \notin (e_{0},\ldots,e_{n},e)$, it is a
nontrivial cycle.  Notice that in the definition of monotone cycle,
the final edge is distinguished.  In particular, a cyclic reordering
of the edges in a monotone cycle may not be a monotone cycle.
\end{remark}

\begin{lemma}\label{lm:mc->asc}
If $\Gamma$ has a strict monotone cycle, then $G$ is ascending.
Further, if $\Gamma$ has a strict monotone cycle, then $\Gamma$ has an
immersed strict monotone cycle.
\end{lemma}

\begin{proof}
If $\Gamma$ contains a strict monotone cycle which is a single edge
$e$, then either $e$ is a strict ascending or strict virtual ascending
loop.  Therefore, after an $\As$--move in the second case, we see
that $G$ is ascending.  Otherwise, suppose $(e_{0},\ldots,e_{n},e)$ is
a strict monotone cycle in $\Gamma$.  Then we can slide $\bar{e}$
along $(e_{0},\ldots,e_{n})$, turning $e$ into a loop.  After the slide
move, the modulus of the loop is a nontrivial integer, hence $e$ is
either a strict ascending or strict virtual ascending loop.  As
before, this shows that $G$ is ascending.

For the second statement in the lemma, we must show that after
tightening, a monotone cycle is still a monotone cycle.  This is clear
since if $(e_{0},\ldots,e_{n},e)$ is a monotone cycle, then $e_{i}
\neq e,\bar{e}$.  Therefore, after tightening, the edge $e$
remains in the cycle and the only tightening occurs in the edge path
$(e_{0},\ldots,e_{n})$, which remains an $\bar{e}$--edge path
after tightening.
\end{proof}

\begin{example}\label{ex:counterexample}
The converse to the first statement of \fullref{lm:mc->asc} does not
hold in general, though we shall prove it in a special case in
\fullref{prop:mc=asc}.  A counterexample is given by the labeled
graphs in \fullref{fig:counterexample}.  The labeled graphs in this
figure represent the same GBS group; the labeled graph on the right is
obtained by sliding $e_{3}$ over $\bar{e}$.  The labeled graph on 
the left contains a strict monotone cycle, namely the cycle
$(e_{1},e_{2},e_{3},e)$.  After sliding $\bar{e}$ over
$(e_{1},e_{2},e_{3})$, $e$ is a virtual ascending loop with labels
$\lambda(e) = 6,\lambda(\bar{e}) = 132$. 

\begin{figure}[ht!]
\labellist
\hair 2pt
\pinlabel {\scriptsize $22$} [br] at 1 30
\pinlabel {\scriptsize $6$} [tl] at 5 26
\pinlabel {\scriptsize $60$} [tr] at 50 26
\pinlabel {\scriptsize $60$} [br] at 49 30
\pinlabel {\scriptsize $6$} [bl] at 55 29
\pinlabel {\scriptsize $15$} [tl] at 55 26
\pinlabel {\scriptsize $12$} [bl] at 105 29
\pinlabel {\scriptsize $5$} [tl] at 105 26
\pinlabel {\small slide} [b] at 158 30
\pinlabel {\scriptsize $6$} [bl] at 244 30
\pinlabel {\scriptsize $60$} [br] at 290 30
\pinlabel {\scriptsize $6$} [br] at 239 29
\pinlabel {\scriptsize $22$} [tr] at 239 26
\pinlabel {\scriptsize $6$} [bl] at 295 29
\pinlabel {\scriptsize $15$} [tl] at 295 26
\pinlabel {\scriptsize $12$} [bl] at 345 29
\pinlabel {\scriptsize $5$} [tl] at 345 26
\pinlabel {$e$} [t] at 27 24
\pinlabel {$e_3$} [b] at 27 56
\pinlabel {$e_2$} [t] at 78 0
\pinlabel {$e_1$} [b] at 78 56
\endlabellist
\begin{center}
\includegraphics{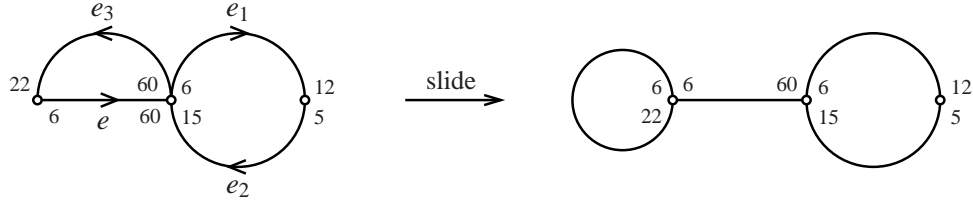}
\end{center}
\caption{Sliding $e_3$ over $\bar{e}$ results in a graph with no
strict monotone cycles.} 
\label{fig:counterexample}
\end{figure}
We claim that the labeled graph on the right does not have any strict
monotone cycles.  First, notice that none of the edges
$e_{1},\bar{e}_1,e_{2},\bar{e}_{2}$ can slide.  Also, since $e$ is
separating (a fact that remains true after sliding $e$ or $\bar{e}$),  no
strict monotone cycle can end with $e$ or $\bar{e}$. Finally, notice that
$\bar{e}_{3}$ cannot slide.  Hence, if there is a strict monotone cycle,
it must be of the form $(\alpha,\bar{e}_{3})$, where $\alpha$ is an
$e_{3}$--edge path.  In particular, $\lambda(e_{3})q(\alpha)$ must be
divisible by $\lambda(\bar{e}_{3}) = 22$.  However, the only place the
prime number $11$ appears in the labeled graph is in the label
$\lambda(\bar{e}_{3})$, and since $\bar{e}_{3} \notin \alpha$,
$\lambda(e_{3})q(\alpha)$ is not divisible by $11$ for any $e_{3}$--edge
path $\alpha$.  Therefore, the labeled graph on the right cannot contain
a strict monotone cycle. 
\end{example}

\begin{remark}
In general, finding a monotone cycle requires a solution to the
conjugacy problem for GBS groups (see \fullref{lm:mobile}). This
problem is not yet known to be solvable.
\end{remark}

\begin{definition}\label{def:slidespace}
Given $\Gamma \in \LG(G)$ and $e \in E(\Gamma)$, we denote by
$\Sl(\Gamma,e) \subseteq \LG(G)$ the set of reduced labeled graphs
obtained from $\Gamma$ by a sequence of slides of $e$ and
$\bar{e}$. $\Sl(\Gamma,e)$ is then called the \emph{slide space of
$e$ (based at $\Gamma$)}. 
\end{definition}

\begin{proposition}\label{prop:infiniteslide}
Let $\Gamma \in \LG(G)$ and $e \in E(\Gamma)$.  Then
$|\Sl(\Gamma,e)| = \infty$ if and only if $\Gamma$ contains a strict
$e$--integer cycle or a strict $\bar{e}$--integer cycle.
\end{proposition}

\begin{proof}
By \fullref{rm:e-edgepath}, it is clear that if $\Gamma$ contains a
strict $e$--integer cycle or a strict $\bar{e}$--integer cycle,
then $|\Sl(\Gamma,e)| = \infty$.

For the converse let $\Gamma_{i}$ be an infinite sequence of labeled
graphs in $\Sl(\Gamma,e)$.  As the number of edges in the graphs
$\Gamma_{i}$ is constant, there is a subsequence such that $\Gamma_{i}
= \Gamma'$ (as \emph{unlabeled graphs}) for some fixed graph
$\Gamma'$.  Thus either $|\lambda_{i}(e)|$ or $|\lambda_{i}(\bar{e})|$
is an unbounded sequence of natural numbers.  By interchanging $e$ for
$\bar{e}$ and passing to a subsequence if necessary, we can
assume that $|\lambda_{i}(e)|$ is a strictly increasing sequence of
natural numbers.  Since slides of $e$ commute with slides of
$\bar{e}$, we can assume that the labeled graphs $\Gamma_{i}$ are
obtained from each other without sliding $\bar{e}$.  There is a
finite number of 
primes appearing in the sequence $\{|\lambda_{i}(e)|\}$.  Indeed, this
list is contained in the set of primes that appear on any labeled
graph for $G$.  Therefore, by the following lemma, there are $n,n'$
such that $\lambda_{n}(e)$ divides $\lambda_{n'}(e)$.  Let $\gamma$ be
the $e$--edge path that $e$ slid along transforming $\Gamma$ into
$\Gamma_{n}$ and $\gamma'$ the strict $e$--integer cycle that $e$ slid
along transforming $\Gamma_{n}$ into $\Gamma_{n'}$.  Then clearly
$\gamma\gamma' \bar{\gamma}$ is a strict $e$--integer cycle in
$\Gamma$.
\end{proof}

\begin{lemma}\label{lm:numbertheory}
Let $\{m_{i}\}$ be a strictly increasing sequence of natural numbers 
such that only finitely many primes appear in the sequence.  Then 
there are distinct indices $n,n'$ such that $m_{n}$ divides $m_{n'}$.   
\end{lemma}

\begin{proof}
We will prove this by induction on the number of primes appearing in 
the sequence $\{m_{i}\}$.  If there is only one prime appearing, then 
the lemma is obvious.  

Suppose that $N$ primes appear in the sequence $\{m_{i}\}$.  To 
any element $m$ in the sequence we associate a point in $\Z^{N}_{\geq 
0}$ (i.e., the first orthant of $\Z^{N}$) by:
\begin{equation*}\label{eq:numbertheory}
    \prod_{j=1}^{N} p_{j}^{k_{j}} \mapsto (k_{1},\ldots,k_{N})
\end{equation*}
where $\prod_{j=1}^{N}p_{j}^{k_{j}}$ is the prime decomposition of
$m$.  For any element $m_i$ in the sequence, we denote the $j$th 
coordinate in this assignment by $(m_i)_{j}$.  If there is some element 
$m_{i}$ such that $(m_{1})_{j} \leq (m_{i})_{j}$ for all $j$,  
then $m_{1}$ divides $m_{i}$ and the conclusion of the lemma
holds. 

Otherwise, by passing to a subsequence, we can assume that 
$(m_{i})_{j} < (m_{1})_{j}$ for some fixed $j$ and all $i$.  By 
further passing to a subsequence we can assume that $(m_{i})_{j} = 
M$ for all $i$.  Then $\{m_{i}/p_{j}^{M}\}$ is a strictly increasing
sequence of natural numbers in which only $N-1$ primes appear.  Now apply
induction to complete the proof. 
\end{proof}

\begin{definition}\label{def:mobile} 
Let $\Gamma \in \LG(G)$.  An edge $e \in E(\Gamma)$ is \emph{mobile}
if either:
\begin{enumerate}
    \item  there is a strict monotone cycle of the form $(e_0,
      \ldots, e_n, e)$ or $(e_0, \ldots, e_n, \bar{e})$; or
    
\item  $|\Sl(\Gamma,e)| = \infty$ (equivalently, by
  \fullref{prop:infiniteslide}, $\Gamma$ contains a strict
  $e$--integer cycle or a strict $\bar{e}$--integer cycle).
\end{enumerate}
An edge that is not mobile is called \emph{non-mobile}. Note that
mobility is a property of geometric edges: $e$ is
mobile if and only $\bar{e}$ is. 
\end{definition}

\begin{remark}\label{rm:mobile}
By \fullref{prop:infiniteslide} there is an algorithm to determine
whether a given edge $e \in E(\Gamma)$ is mobile or not.  Indeed,
given an edge we can start making an exhaustive search of
$\Sl(\Gamma,e)$.  Either this space is finite or we find an strict
$e$--integer cycle or strict $\bar{e}$--integer cycle.  In the latter
case, $e$ is mobile.  If the slide space is finite, we can search
these graphs to see if $e$ is a strict ascending or strict virtual
ascending loop in any of the graphs.  An affirmative answer implies
that $e$ is mobile, a negative answer implies that $e$ is non-mobile.
\end{remark}

Let $T$ denote the Bass--Serre tree covering $\Gamma$.

\begin{lemma}\label{lm:mobile}
An edge $e \in E(\Gamma)$ is mobile if and only if $G_{\tilde{e}}^{t}
\subsetneq G_{\tilde{e}}$ for some $t \in G$ and some lift $\tilde{e} \in
E(T)$ of $e$. 
\end{lemma}

\begin{proof}
It is clear that if $e$ or $\bar{e}$ is the last edge of a strict
monotone cycle, or there is strict $e$--integer cycle or strict
$\bar{e}$--integer cycle then there is a lift $\tilde{e}$ and a
$t \in G$ satisfying the conclusion of the lemma.

For the converse, given $t$ and $\tilde{e}$ with $G_{\tilde{e}}^{t}
\subsetneq G_{\tilde{e}}$, we can replace $t$ if needed to arrange that there
are no $G$--translates of $\tilde{e}$ along the edge path connecting
$\tilde{e}$ to $t\tilde{e}$. Also we can assume without
loss of generality that $\iv(\tilde{e})$ separates $\tilde{e}$ from
$t\tilde{e}$. Let $\tilde{\alpha}$ be the path in $T$ from $\tilde{e}$ 
to $t \tilde{e}$, and let $\alpha$ be its image in $\Gamma$. 

If $\tilde{e}$ and $t\tilde{e}$ are coherently oriented, then
$(\alpha, \bar{e})$ is a strict monotone cycle.  Otherwise,
$\alpha$ is a strict $e$--integer cycle.
\end{proof}

\begin{corollary}\label{co:mobile}
If $\Gamma, \Gamma' \in \LG(G)$ are related by slide moves and $e
\in E(\Gamma)$ is mobile, then $e$ is mobile in $\Gamma'$. 
\end{corollary}

\begin{proof}
This follows from \fullref{lm:mobile} since edge stabilizers are
unchanged by slide moves. 
\end{proof}

\fullref{ex:counterexample} shows that both parts of the
definition for mobility are needed.  The edge $e$ is part of a strict
monotone cycle $(e_{1},e_{2},e_{3},e)$ in the labeled graph on the
left and hence is mobile.  In the labeled graph on the right, as noted in
the example, there are no strict monotone cycles, but there is a
strict $\bar{e}$--integer cycle, so $e$ is mobile.

\begin{remark}\label{rm:nonmobilecanonical}
The set of non-mobile edges is preserved by slides, inductions, and
$\As^{\pm 1}$--moves.  To make sense of the third case, observe that
even though an $\As^{\pm 1}$--move changes the set of edges, the
edges directly involved in the move are all mobile, so each non-mobile
edge is present before and after the move, and its status does not
change (by \fullref{lm:mobile}).  In the case of an induction move,
the loop is mobile before and after, and mobility of other edges is
not affected, again by \fullref{lm:mobile}.  Therefore, for any
labeled graph space, we can compare non-mobile edges between any two
labeled graphs. 
\end{remark}

\begin{lemma}\label{lm:nonmobileovermobile} 
In a labeled graph, a non-mobile edge cannot slide over a mobile edge.
\end{lemma}

\begin{proof}
Suppose an edge $f$ slides over a mobile edge $e$.  Then there are
lifts $\tilde{f}$ and $\tilde{e}$ in the covering tree $T$ such that
$\tilde{f}$ slides over $\tilde{e}$, and so $G_{\tilde{f}} \subseteq
G_{\tilde{e}}$.  Let $n = [G_{\tilde{e}} : G_{\tilde{f}}]$.  By
\fullref{lm:mobile} there is a $t\in G$ such that $G^t_{\tilde{e}}
\subsetneq G_{\tilde{e}}$.  Let $m = [G_{\tilde{e}} :
G^t_{\tilde{e}}]$.  Then $G^t_{\tilde{f}} \subseteq G^t_{\tilde{e}}
\subsetneq G_{\tilde{e}}$, and $G^t_{\tilde{f}}$ is the unique
subgroup of $G_{\tilde{e}}$ of index $mn$.  This implies that
$G^t_{\tilde{f}}$ is the subgroup of $G_{\tilde{f}}$ of index $m$,
which is greater than $1$, and so $f$ is mobile.
\end{proof}

\subsection{Slide relations}\label{ssc:relations}

In this subsection we will work out some methods to rearrange sequences
of slide moves.  In particular, we will show that any sequence of
slides can be rewritten so that non-mobile edges slide before mobile
edges, and individual non-mobile edges can be slid one at a time.  To
simplify the discussion, we will only consider positive labeled graphs. 
All slides in this section are between \emph{reduced} trees (that is, the
slides take place ``in $\LG^+(G)$'').  

\begin{notation}\label{no:slides}
If $\Gamma \in \LG(G)$, $e \in E(\Gamma)$ and $A$ is an $e$--edge
path, we will use the notation $e/A$ to denote the slide move of $e$
over $A$.  When we write a composition of slides $e/A
\cdot f/B$ we will always assume that $f/B$ is a valid slide move after
sliding $e$ over $A$.  We have some obvious relations: $e/A \cdot e/A' =
e/AA'$ and $e/\bar{A}$ is the inverse of $e/A$ (here $\bar{A}$
is the reverse of the path $A$). 

Throughout the rest of the section, $A$ denotes an $e$ or
$\bar{e}$--edge path and $B$ denotes an $f$ or $\bar{f}$--edge
path. Likewise for $A'$, $B'$, etc.  We will use $\alpha$ to denote an 
$e$ or $\bar{e}$--edge path not containing $f$ or $\bar{f}$, and $\beta$
an $f$ or $\bar{f}$--edge path not containing $e$ or $\bar{e}$. 
\end{notation}

The following proposition is our current goal.

\begin{proposition}\label{prop:commute}
Suppose $\Gamma \in \LG^+(G)$ and $e,f \in E(\Gamma)$ $(e \neq
f,\bar{f})$ where $f$ is non-mobile.  Suppose $e/A \cdot f/B $ is
valid slide sequence in $\LG^+(G)$. Then:
\begin{equation*}
    e/A \cdot f/B \ = \ f/B' \cdot \bar{f}/B'' \cdot e/A' \cdot
    \bar{e}/A''
\end{equation*}
for some appropriate edge paths $B',B'',A'$ and $A''$.
\end{proposition}

We will establish this proposition by a careful analysis of how to
commute individual slide moves past one another.  We begin by listing
several basic relations. 

\begin{definition}
In some of the slide relations below, \emph{renaming} occurs. This does
not mean that the edges themselves are renamed. Rather, when the
relation is used to substitute some slide moves for others inside a
larger sequence of moves, the moves in the larger sequence occurring
\emph{after} the newly substituted moves need to be renamed, so that they
still refer to the same edges as before. For example, the instruction
``rename $e \mapsto f$, $f \mapsto \bar{e}\,$'' means that moves such
as $e/\alpha$, $f/\bar{e}$, $e/\bar{f}$ occurring later in the
sequence should now be written as $f/\alpha$,
$\bar{e}/\bar{f}$, $f/e$. The reason for this should become
clear in the proof of the next lemma. 
\end{definition}

\begin{lemma}\label{lm:relations}
Suppose that $f$ is non-mobile.  Then the following relations are 
valid:
\begin{enumerate}
\item    $e/\alpha \cdot f/\beta \ = \ f/\beta \cdot e/\alpha$\label{al:ab} 
\item    $e/\alpha \cdot f/e \ = \ f/\bar{\alpha}e \cdot 
    e/\alpha$\label{al:ae} 
\item    $e/\alpha \cdot f/\bar{e} \ = \ f/\bar{e}\alpha \cdot 
    e/\alpha$\label{al:aE} 
\item    $e/f \cdot f/\beta \ = \ f/\beta \cdot e/\beta f$\label{al:fb} 
\item    $e/f \cdot f/e \ = \ f/e$\label{al:fe} 
\item    $e/f \cdot f/\bar{e} \ = \ e/f \cdot \bar{e}/f, \mbox{ then rename } 
    e\mapsto \bar{f},f \mapsto \bar{e}$\label{al:fE} 
\item    $e/\bar{f} \cdot f/\beta \ = \ f/\beta \cdot 
    e/\bar{f}\bar{\beta}$\label{al:Fb} 
\item    $e/\bar{f} \cdot f/e \ = \ \bar{f}/e, \mbox{ then rename }
    e\mapsto f, f \mapsto \bar{e}$\label{al:Fe}
\item    $e/\bar{f} \cdot f/\bar{e} \ = \ e/\bar{f}$\label{al:FE}
\end{enumerate}
where $f,\bar{f} \notin \alpha$ and $e,\bar{e} \notin
\beta$.  Furthermore, after substituting and renaming moves, $f$ still
refers to a non-mobile edge. 
\end{lemma}

\begin{proof}
In the diagrams below, the heavy edge is $f$ and the light edge is
$e$. Note that in cases \eqref{al:fE} and \eqref{al:Fe}, later references
to these edges will be renamed. Since $f$ slides over $e$ or
$\bar{e}$ in these cases, $e$ must be non-mobile by
\fullref{lm:nonmobileovermobile}. So it remains true that $f$ is
non-mobile in later moves, after the renaming step. 

Now consider the individual cases, recalling that $\alpha$ and $\beta$ do
not contain $e$, $\bar{e}$, $f$, or $\bar{f}$. Case
\eqref{al:ab} is obvious. Case \eqref{al:ae} is clear after noting that
$\iv(e) = \iv(\alpha)$ and $\tv(\alpha) = \iv(f)$. In case \eqref{al:aE}
we have $\iv(\alpha) = 
\iv(e)$ and $\tv(e) = \iv(f)$ and the relation is clear. In case
\eqref{al:fb} we have $\iv(e) = \iv(f) = \iv(\beta)$ and the relation is
clear. 

For \eqref{al:fe}, shown below, $f$ is a loop at $\iv(e)$. The labels
$\lambda(e), \lambda(f)$ are of the form $ca, a$ since $e$ slides over
$f$. Then, since $f$ slides over $e$, we must have $cb \mid a$ (where $b
= \lambda(\bar{f})$), hence $b \mid a$. Since $f$ is non-mobile, we
then have $b=a$. Hence the first slide may simply be omitted. 
\begin{center}
\labellist
\small \hair 2pt
\pinlabel {$ca$} [bl] at 22 43
\pinlabel {$a$} [br] at 18 28
\pinlabel {$b$} [bl] at 22 28
\pinlabel {$cb$} [bl] at 135 43
\pinlabel {$a$} [br] at 131 28
\pinlabel {$b$} [bl] at 135 28
\pinlabel {$cb$} [bl] at 261 17
\pinlabel {$b$} [tr] at 257 14
\endlabellist
\includegraphics{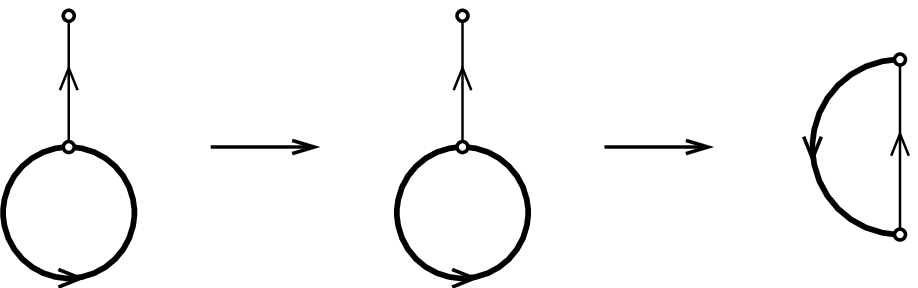}
\end{center}

For \eqref{al:fE}, shown below, $e$ is a loop at $\iv(f)$ and we have
$\lambda(e), \lambda(f)$ of the form $ca, a$ as before. Since $f$ slides
over $\bar{e}$, we have $b \mid a$ (where $b =
\lambda(\bar{e})$) and the new label of $f$ becomes $acd/b$ (where
$d = \lambda(\bar{f})$). This integer is divisible by $d$, and so
$f$ is now virtually ascending. Since $f$ is non-mobile, we conclude that
$acd/b = d$, so $ac  = b$. Since $b \mid a$, we now have $c = 1$ and
$a=b$. The result of the two moves can now be achieved by sliding $e$ and
$\bar{e}$ over $f$. After this move, $\bar{f}$ is in the position
previously occupied by $e$, so later references to $e$ should be renamed
as $\bar{f}$. Similarly, references to $f$ should be renamed to
$\bar{e}$ ($e$ would work equally well in this case). 
\begin{center}
\labellist
\small \hair 2pt
\pinlabel {$a$} [bl] at 21 17
\pinlabel {$d$} [br] at 68 17
\pinlabel {$a$} [bl] at 141 17
\pinlabel {$d$} [br] at 187 17
\pinlabel {$d$} [br] at 307 27
\pinlabel {$d$} [bl] at 310 27
\pinlabel {$ca$} [br] at 20 31
\pinlabel {$b$} [bl] at 22 31
\pinlabel {$b$} [br] at 138 30
\pinlabel {$cd$} [bl] at 191 30
\pinlabel {$b$} [br] at 258 43
\pinlabel {$cd$} [bl] at 311 43
\endlabellist
\includegraphics{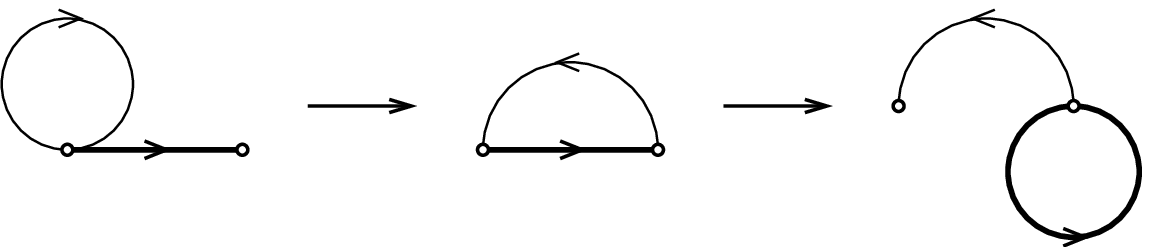}
\end{center}
In case \eqref{al:Fb} we have $\iv(e) = \tv(f)$ and $\iv(f) = \iv(\beta)$
and the relation is 
clear. For \eqref{al:Fe} we have $\iv(e) = \tv(f)$, and $\lambda(e),
\lambda(\bar{f})$ of the form $ba, a$. After the first slide
$\lambda(e)$ becomes $bc$ where $c = \lambda(f)$, and since the second
slide occurs we have that $bc \mid c$. Hence $b=1$ and we originally have
$\lambda(e) = \lambda(\bar{f})$. Now the same labeled graph results
by sliding $\bar{f}$ over $e$. In later moves, $e$ should be renamed
as $f$, and $f$ as $\bar{e}$, since $f$ and $\bar{e}$ now
occupy the previous positions of $e$ and $f$.  

Case \eqref{al:FE} is shown below: 
\begin{center}
\labellist
\small \hair 2pt
\pinlabel {$a$} [bl] at 6 31
\pinlabel {$b$} [br] at 49 31
\pinlabel {$ca$} [bl] at 6 42
\pinlabel {$d$} [br] at 49 42
\pinlabel {$a$} [bl] at 120 17
\pinlabel {$b$} [bl] at 168 17
\pinlabel {$cb$} [br] at 166 32
\pinlabel {$d$} [bl] at 169 32
\pinlabel {$a$} [bl] at 252 17
\pinlabel {$b$} [bl] at 300 17
\endlabellist
\includegraphics{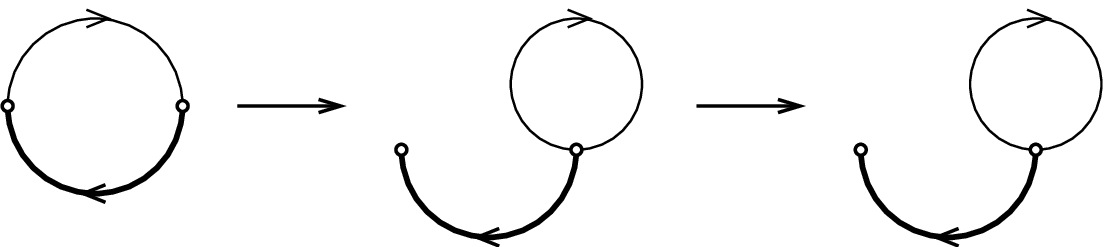}
\end{center}
The labels $\lambda(e), \lambda(\bar{f})$ are of the form $ca, a$
and $\lambda(e)$ becomes $cb$ after the first slide (where $b =
\lambda(f)$). From the second slide we deduce that $d =
\lambda(\bar{e})$ divides $b$. Now $e$ is virtually ascending, and
it is non-mobile since $f$ can slide over it. So $cb = d$ and the second
slide may be omitted. 
\end{proof}

The next result is a straightforward application of the relations
\eqref{al:ab}--\eqref{al:fb} of \fullref{lm:relations}. It will be used
to establish a special case of \fullref{prop:commute}, when 
either $f,\bar{f} \notin A$ or $e,\bar{e} \notin B$. 

\begin{lemma}\label{lm:easycase}
Suppose that $f$ is non-mobile. Then the following relations are valid: 
\begin{enumerate}
\item $e/\alpha f \alpha' \cdot f/\beta \ = \ f/\beta \cdot e/\alpha
  \beta f \alpha'$ \label{al:easy1} 
\item $e/\alpha \bar{f} \alpha' \cdot f/\beta \ = \ f/\beta \cdot
  e/\alpha \bar{f} \bar{\beta} \alpha'$ \label{al:easy2} 
\item $e/\alpha \cdot f/\beta e \beta' \ = \ f/\beta
  \bar{\alpha}e\beta' \cdot e/\alpha$ \label{al:easy3} 
\item $e/\alpha \cdot f/\beta \bar{e} \beta' \ = \
  f/\beta\bar{e}\alpha\beta' \cdot e/\alpha$ \label{al:easy4} 
\end{enumerate}
\end{lemma}

\begin{proof}
For \eqref{al:easy1} we write 
\begin{align*}
e/\alpha f \alpha' \cdot f/\beta \ &= \ e/\alpha f \cdot f/\beta \cdot 
e/\alpha' \tag*{\scriptsize \quad By \ref{lm:relations}(\ref{al:ab}):}\\
&= \ e/\alpha \cdot f/\beta \cdot e/\beta f \alpha' \tag*{\scriptsize
\quad \ref{lm:relations}(\ref{al:fb}):}\\ 
&= \ f/\beta \cdot e/\alpha \beta f \alpha'. \tag*{\scriptsize \quad
\ref{lm:relations}(\ref{al:ab}):} 
\end{align*}
The other relations are similar. 
\end{proof}

The next relations will form the basis of the proof of
\fullref{prop:commute}. 

\begin{lemma}\label{lm:hardcase}
Suppose that $f$ is non-mobile.  Then the following relations are 
valid:     
\begin{enumerate}
\item $e/\alpha f \alpha' \cdot f/\beta e \beta' \ = \ f/\bar{\alpha}e \beta' 
\cdot e/\alpha\beta$\label{al:afabeb} 
\item $e/\alpha f \alpha' \cdot f/\beta \bar{e}\beta' \ = \ f/\beta \cdot
\bar{f}/\alpha' \cdot e/\alpha\beta f\bar{\alpha}' \cdot
\bar{e}/f \beta', \mbox{ rename } e \mapsto \bar{f},f
\mapsto \bar{e}$\label{al:afabEb} 
\item $e/\alpha \bar{f} \alpha' \cdot f/\beta \bar{e} \beta' \ =
 \ f/\alpha' \beta' \cdot e/\alpha\bar{f}\bar{\beta}'$\label{al:aFabEb} 
\item $e/\alpha \bar{f} \alpha' \cdot f/\beta e \beta' \ = \ f/\alpha'
\cdot \bar{f}/\bar{\alpha}e \cdot e/\alpha \cdot
\bar{e}/\beta', \mbox{ rename } e \mapsto f,f \mapsto
\bar{e}$\label{al:aFabeb} 
\end{enumerate}
\end{lemma}

\begin{proof}
The first three of these follow from straightforward computations,
similar to \fullref{lm:easycase}. 
\begin{align}
e/\alpha f \alpha' \cdot f/\beta e \beta' \ &= \ e/\alpha \cdot f/\beta 
\cdot e/\beta f \cdot f/\bar{\alpha}'e \cdot e/\alpha' \cdot
f/\beta' \tag*{\scriptsize \quad By
\ref{lm:relations}(\ref{al:ab},\ref{al:fb},\ref{al:ae}):} \\ 
&= \ f/\beta \cdot e/\alpha\beta \cdot f/\bar{\alpha}' \cdot
e/\bar{\alpha}'f \cdot f/e\beta' \cdot e/\alpha' \tag*{\scriptsize \quad 
  \ref{lm:relations}(\ref{al:ab},\ref{al:fb}):} \\
&= \ f/\beta\bar{\alpha}' \cdot e/\alpha\beta \cdot f/\alpha'e \cdot
e/\bar{\alpha}' \cdot f/\beta' \cdot e/\alpha' 
\tag*{\scriptsize \quad 
\ref{lm:relations}(\ref{al:ab},\ref{al:fe},\ref{al:ae}):} \\ 
&= \ f/\beta \cdot e/\alpha\beta \cdot f/e\beta' 
\tag*{\scriptsize \quad \ref{lm:relations}(\ref{al:ab}) and
cancellation:}\\  
&= \ f/\bar{\alpha}e\beta' \cdot e/\alpha\beta.
\tag*{\scriptsize \quad \ref{lm:relations}(\ref{al:ae},\ref{al:ab}):} 
\end{align}
This proves \eqref{al:afabeb}. For \eqref{al:afabEb} we have:  
\begin{align}
e/\alpha f\alpha' \cdot f/\beta\bar{e}\beta' \ &= \ e/\alpha f \cdot 
f/\beta \cdot e/\alpha' \cdot f/\bar{e}\beta' \notag \\
&= \ f/\beta \cdot e/\alpha\beta f \cdot f/\bar{e}\alpha'\beta' \cdot 
e/\alpha'\tag*{\scriptsize \quad 
\ref{lm:relations}(\ref{al:fb},\ref{al:aE},\ref{al:ab}):} \\
&= \ f/\beta \cdot e/\alpha\beta f \cdot \bar{e}/f \cdot 
\bar{e}/\alpha'\beta' \cdot \bar{f}/\alpha', \text{ rename } e \mapsto
\bar{f}, f \mapsto \bar{e} \tag*{\scriptsize \quad 
\ref{lm:relations}(\ref{al:fE}):} \\
&= \ f/\beta \cdot e/\alpha\beta f \cdot \bar{e}/f \cdot 
\bar{f}/\alpha' \cdot \bar{e}/\alpha'\beta', \text{ rename } e \mapsto
\bar{f}, f \mapsto \bar{e}  \tag*{\scriptsize \quad 
\ref{lm:relations}(\ref{al:ab}):} \\
&= \ f/\beta \cdot e/\alpha\beta f \cdot \bar{f}/\alpha' \cdot 
\bar{e}/f \beta', \text{ rename } e \mapsto \bar{f}, f
\mapsto \bar{e}  \tag*{\scriptsize \quad 
\ref{lm:relations}(\ref{al:Fb}) and cancellation:} \\
&= \ f/\beta \cdot \bar{f}/\alpha' \cdot e/\alpha\beta f \bar{\alpha}'
\cdot \bar{e}/f\beta', \text{ rename } e \mapsto \bar{f}, f
\mapsto \bar{e}. \tag*{\scriptsize \quad 
\ref{lm:relations}(\ref{al:Fb},\ref{al:ab}):}
\end{align}
Note that part of the third line has undergone renaming. The renaming
instruction is still needed for any subsequent moves. 
Next consider \eqref{al:aFabEb}: 
\begin{align}
e/\alpha\bar{f}\alpha' \cdot f/\beta\bar{e}\beta' \
&= \ e/\alpha\bar{f} \cdot f/\beta \cdot e/\alpha' \cdot f/\bar{e}\beta'
\notag \\ 
&= \ f/\beta \cdot
e/\alpha\bar{f}\bar{\beta}\alpha' \cdot f/\bar{e}\beta' \tag*{\scriptsize
\quad \ref{lm:relations}(\ref{al:Fb},\ref{al:ab}):}\\
&= \ f/\beta \cdot e/\alpha\bar{f}\bar{\beta} \cdot 
f/\bar{e}\alpha'\beta' \cdot e/\alpha'
\tag*{\scriptsize \quad \ref{lm:relations}(\ref{al:aE},\ref{al:ab}):}\\ 
&= \ f/\beta \cdot e/\alpha\bar{f} \cdot 
f/\bar{e}\bar{\beta}\alpha'\beta' \cdot e/\bar{\beta}\alpha'
\tag*{\scriptsize \quad \ref{lm:relations}(\ref{al:aE},\ref{al:ab}):}\\ 
&= \ f/\beta \cdot e/\alpha \cdot f/\bar{\beta} \cdot e/\bar{f} \beta
\cdot f/\alpha'\beta' \cdot e/\bar{\beta}\alpha' 
\tag*{\scriptsize \quad \ref{lm:relations}(\ref{al:FE},\ref{al:Fb}):}\\ 
&= \ e/\alpha\bar{f} \cdot f/\alpha'\beta' \cdot e/\alpha \tag*{\scriptsize
\quad \ref{lm:relations}(\ref{al:ab}) and cancellation:} \\
&= \ f/\alpha'\beta' \cdot e/\alpha\bar{f}\bar{\beta}'. \tag*{\scriptsize
\quad \ref{lm:relations}(\ref{al:Fb},\ref{al:ab}) and cancellation:}
\end{align}

Finally we prove \eqref{al:aFabeb}.  Notice that as
$\alpha\bar{f}\alpha'$ is an $e$--edge path $\iv(\alpha') = \iv(f)$
and after sliding $e$ we have that $\iv(e) = \tv(\alpha')$.  Also,
since $\beta e\beta'$ is an $f$--edge path after sliding $e$,
$\iv(\beta) = \iv(f) = \iv(\alpha')$ and $\tv(\beta) = \iv(e) =
\tv(\alpha')$ Therefore, as neither $\alpha'$ nor $\beta$ contain
$e,\bar{e}$, $\beta\bar{\alpha}'$ is a cycle before sliding $e$.
Since after sliding $e$ over $\alpha\bar{f}\alpha'$ we can slide $f$
over $\beta e$ we have that
$\lambda_{\Gamma}(e)q_{\Gamma}(\alpha,\bar{f},\alpha')$ divides
$\lambda_{\Gamma}(f)q_{\Gamma}(\beta)$. (Here $\Gamma$ is the labeled
graph just before the slide moves under discussion.)  In particular, 
after sliding $f$
along $\beta$, we can slide it back along $\bar{\alpha}'$.  Finally,
since $\lambda_{\Gamma}(f)$ divides
$\lambda_{\Gamma}(e)q_{\Gamma}(\alpha,\bar{f})$ which divides
$\lambda_{\Gamma}(f)q_{\Gamma}(\beta\bar{\alpha}')$, we have that 
$q_{\Gamma}(\beta\bar{\alpha}')$ is an integer.  As $f$ is non-mobile,
this integer must be $1$ (recall that we are assuming that all labels
are positive).  Hence $f/\beta = f/\alpha'$.  Now it is easy to verify
that \eqref{al:aFabeb} is a valid relation: 
\begin{align}
e/\alpha\bar{f}\alpha' \cdot f/\alpha'e\beta' \ &= \ e/\alpha \bar{f} 
\cdot f/\alpha' \cdot e/\alpha' \cdot f/e\beta' \notag \\ 
&= \ f/\alpha' \cdot e/\alpha\bar{f} \cdot f/e\beta' \tag*{\scriptsize
\quad \ref{lm:relations}(\ref{al:Fb},\ref{al:ab}) and cancellation:} \\ 
&= \ f/\alpha' \cdot e/\alpha \cdot \bar{f}/e \cdot 
\bar{e}/\beta', \text{ rename } e \mapsto {f},f \mapsto \bar{e}
\tag*{\scriptsize
\quad \ref{lm:relations}(\ref{al:Fe}):} \\ 
&= \ f/\alpha' \cdot \bar{f}/\bar{\alpha}e \cdot e/\alpha \cdot 
\bar{e}/\beta', \text{ rename } e \mapsto {f},f \mapsto
\bar{e}. \tag*{\scriptsize 
\quad \ref{lm:relations}(\ref{al:ae}):} 
\end{align}
This completes the proof. 
\end{proof}

We are now in a position to prove \fullref{prop:commute}. 

\begin{proof}[Proof of \fullref{prop:commute}]
To simplify the discussion we introduce a shorthand for slide
sequences. Slides of the form $e/\alpha$ or $\bar{e}/\alpha$ are denoted
by $E$, and those of the form $e/\alpha f \alpha'$, $e/\alpha \bar{f}
\alpha'$, $\bar{e}/\alpha f \alpha'$ or $\bar{e}/\alpha \bar{f} \alpha'$
by $E_F$. Likewise define the symbols $F$ and $F_{E}$. Given a slide
sequence, let $m$ denote the number of slides of the form $E_F$ or
$F_E$. Let $n$ denote the number of transitions of the form $E_F F_E$
after omitting the symbols $E,F$. The \emph{complexity} of the sequence
is the pair $(m,n)$, ordered lexicographically. 

We are given the sequence $e /A \cdot f/ B$, which decomposes into a
slide sequence consisting of $E$'s and $E_F$'s, followed by $F$'s and
$F_E$'s. Our strategy is to apply slide relations to reduce complexity,
until $n=0$. If $n=0$ then we have a sequence in which no $E_F$ appears
before an $F_E$. To complete the argument in this case,
\fullref{lm:relations}\eqref{al:ab} will transform any $EF$ to $FE$;
\fullref{lm:easycase}(\ref{al:easy1},\ref{al:easy2}) transforms any $E_F F$
to $F E_F$; and \fullref{lm:easycase}(\ref{al:easy3},\ref{al:easy4})
transforms any $E F_E$ to $F_E E$. Using these relations, the sequence
can be transformed to one consisting of $F$'s and $F_E$'s followed by
$E$'s and $E_F$'s. Lastly, since slides of $e$ and $\bar{e}$
(respectively, $f$ and $\bar{f}$) commute, the sequence can be put into
the desired form $f/B' \cdot \bar{f}/B'' \cdot e/A' \cdot \bar{e}/A''$. 

Next we show how to reduce complexity if $n > 0$. We will be
applying the relations of \fullref{lm:hardcase}, some of which involve
renaming. When this occurs, the symbols $E_F$ and $F_E$, and the
symbols $E$ and $F$, will be exchanged throughout part of the
sequence. Notice that this in itself does not change $m$. Notice also
that the relations in \ref{lm:hardcase}\eqref{al:afabeb},
\ref{lm:hardcase}\eqref{al:aFabEb} and \ref{lm:hardcase}\eqref{al:aFabeb} 
all reduce $m$.  

There is one additional rewriting move which has not yet been
discussed. The moves $E E_F$ may be rewritten either as $E_F$ or as
$E_F E$, depending on whether the edge $e$ appears with the same
orientation in the two moves. Similarly, $F_E F$ can be rewritten as
$F_E$ or $F F_E$. 

The procedure is first to push all $F$'s to the beginning of the
sequence and all $E$'s to the end, using this last observation and Lemmas 
\ref{lm:relations}\eqref{al:ab} and \ref{lm:easycase}. This does not
change complexity. Then apply relation
\ref{lm:hardcase}\eqref{al:afabeb}, \ref{lm:hardcase}\eqref{al:aFabEb} or
\ref{lm:hardcase}\eqref{al:aFabeb}, if possible, to one of the $E_F F_E$
pairs, to reduce complexity. If none of these apply, then every $E_F F_E$
pair matches the left hand side of relation
\ref{lm:hardcase}\eqref{al:afabEb}. Using this relation does not
obviously reduce complexity, but we can proceed as follows.  

Starting with the rightmost $E_F F_E$ pair, the slide sequence has the
form 
\[ \cdots {(E_F F_E)} {(F_E)}^* {(E_F)}^* {(E)}^* \]
where $^*$ denotes zero or more copies of the symbol. Applying 
\ref{lm:hardcase}\eqref{al:afabEb} to this pair, the sequence becomes 
\[ \cdots {(F F E_F E_F)} {(E_F)}^* {(F_E)}^* {(F)}^*, \]
with no change to the symbols that are not shown. If the
$(E_F)^*$ term in the original sequence is empty then $n$ decreases 
and $m$ stays the same, and complexity has been reduced. Otherwise the new
sequence has the same complexity. If this occurs, apply 
\fullref{lm:hardcase} to 
the newly created rightmost $E_F F_E$ pair. If case \eqref{al:afabeb},
\eqref{al:aFabEb} or \eqref{al:aFabeb} applies, complexity is reduced as
before. If case \eqref{al:afabEb} applies then we are in the situation
just discussed, with empty $(E_F)^*$ term, and $n$ decreases. Thus, in
all cases, complexity has been reduced. 
\end{proof}

The corollary below follows directly, by repeated application of
\fullref{prop:commute}. 

\begin{corollary}\label{co:finitefirst}
Suppose $\Gamma,\Gamma'$ are related by a sequence of slides in $\LG^+(G)$
and $f \in E(\Gamma)$ is non-mobile.  Then there is a labeled graph
$\Gamma_f\in \Sl(\Gamma,f)$ and a sequence of slides $\Gamma_{f} \to
\Gamma'$ in $\LG^+(G)$ during which the edges $f, \bar{f}$ remain
stationary. 
Moreover, if a geometric edge $e, \bar{e} \in E(\Gamma)$ was
stationary in the original slide sequence, then the sequence $\Gamma_f
\to \Gamma'$ may be chosen to leave $e, \bar{e}$ stationary 
as well. 
\end{corollary}

\subsection{Finiteness of $\LG(G)$}\label{ssc:finite}
We can now prove \fullref{th:finite}, along with some
applications. Here is a restatement of the theorem. 

\begin{theorem}
Let $\Gamma \in \LG(G)$, where $G \not= BS(1,n)$.  Then $|\LG(G)| =
\infty$ if and only if $\Gamma$ has a mobile edge.
\end{theorem}

\begin{proof}
Let $G^{+}$ be the GBS group represented by the labeled graph
$(\Gamma,|\lambda|)$. Changing the signs of a labeling has
no effect on divisibility relations, and hence has no effect on slide
moves or mobility of edges. Moreover, the absolute value map $\LG(G) \to
\LG^+(G^{+})$ is finite-to-one, so $|\LG(G)|$ is finite if and only if
$|\LG^+(G^{+})|$ is. Thus, without loss of generality, we may assume that
$\Gamma$ is a positive labeled graph, and we may work in $\LG^+(G)$,
where \fullref{co:finitefirst} is valid. 
    
Suppose $\Gamma$ has a mobile edge $e$.  If there is a strict
monotone cycle then $G$ is ascending, and since $G \not=
BS(1,n)$, it follows that $|\LG(G)| = \infty$.  Otherwise
$|\Sl(\Gamma,e)| = \infty$, which implies that $|\LG(G)| = \infty$.

Next suppose that $\Gamma$ has no mobile edges. In particular, $G$ is
non-ascending.  By \fullref{co:slides}, $\LG^+(G)$ is connected by
slide moves.  Let $e_1, e_2, 
\ldots, e_k$ be the geometric edges of $\Gamma$.  Given any $\Gamma'\in
\LG^+(G)$, \fullref{co:finitefirst} implies that 
there is a sequence of labeled graphs $\Gamma = \Gamma^0,
\Gamma^1, \ldots, \Gamma^k = \Gamma'$ such that $\Gamma^i$ is in the
slide space $\Sl(\Gamma^{i-1},e_i)$ for each $i$.  Since no $\Gamma^i$ has
a mobile edge, these slide spaces are all finite, and therefore
$\LG^+(G)$ is finite. 
\end{proof}

\begin{remark}\label{rm:finite}
Since we have an algorithm to determine whether a given labeled
graph has a mobile edge (\fullref{rm:mobile}), the finiteness
criterion above can be checked algorithmically. 
\end{remark}

\begin{figure}[ht!]
\labellist
\small \hair 2pt
\pinlabel {\scriptsize $2$} [br] at 9 45
\pinlabel {\scriptsize $4$} [bl] at 13 51
\pinlabel {\scriptsize $18$} [tr] at 9 7
\pinlabel {\scriptsize $24$} [t] at 18 0
\pinlabel {\scriptsize $18$} [bl] at 52 27
\pinlabel {\scriptsize $24$} [tl] at 52 24
\endlabellist
\begin{center}
\includegraphics{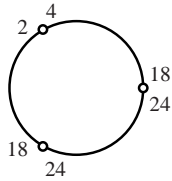}
\end{center}
\caption{A GBS group with finite labeled graph space and integral modulus.} 
\label{fig:nomobiles}
\end{figure}

\begin{example}\label{ex:newexample}
\fullref{fig:nomobiles} shows a labeled graph with modulus a
nontrivial integer. For this reason, the finiteness theorem of Forester
\cite{ar:F06} does not apply. There is only 
one possible slide move, and the only slide afterwards is its
reverse. It follows that there are no mobile edges, by 
\fullref{rm:mobile}. Hence this GBS group has only finitely many reduced
labeled graphs representing it. 
\end{example}

As a consequence of \fullref{th:finite}, we have the following
theorem about the finiteness properties of the group of outer
automorphisms of a GBS group. The proof is exactly as in
Clay \cite[Theorem~75]{ar:CT} or Levitt
\cite[Theorem~1.5]{ar:L07}. Recall that a 
group is \emph{of type $\textup{F}_{\! \infty}$} if it is the fundamental
group of an aspherical cell complex having finitely many cells in each
dimension.  

\begin{theorem}\label{th:fp}
If a GBS group $G$ is represented by a labeled graph that does not
have any mobile edges, then $\Out(G)$ is of type $\textup{F}_{\!
\infty}$. 
\end{theorem}

Another application concerns the isomorphism problem for GBS groups. 

\begin{theorem}\label{th:iso-finite}
There is an algorithm which, given two labeled graphs, one of which
does not have any mobile edges, determines whether the corresponding
GBS groups are isomorphic.
\end{theorem}

\begin{proof}
Let $\Gamma, \Gamma'$ be reduced labeled graphs with corresponding GBS
groups $G, G'$, where $\Gamma$ has no mobile edges. 
\fullref{rm:nonmobilecanonical} implies that the deformation space of
$\Gamma$ is non-ascending. Hence, by \fullref{co:slides}, reduced
trees in this deformation space are connected by slide moves (between
reduced trees). Since $\LG(G)$ is finite (\fullref{th:finite}),
it can be enumerated effectively, by performing all possible slide
sequences, exactly as in the proof of \cite[Corollary 8.3]{ar:F06}. Then 
$G$ and $G'$ are isomorphic if and only if the labeled graph $\Gamma'$ is
found. 
\end{proof}

\section{Betti number one graphs}\label{sc:bettione}

Given a non-elementary GBS group $G$, all labeled graphs $\Gamma\in
\LG(G)$ have the same first Betti number $b(\Gamma)$, since this is
simply the rank of the quotient of $G$ by the subgroup generated by the
elliptic elements. Alternatively, collapse and expansion moves induce 
homotopy equivalences of the underlying graph. Thus we also denote this
number $b(G)$. In this section we will only consider GBS groups $G$ such
that $b(G) = 1$. As before, all GBS groups in this section are assumed to
be non-elementary. 

\begin{remark}\label{rm:monotonecheck}
Suppose $\Gamma \in \LG(G)$ with $b(\Gamma) = 1$.  If there is
a strict monotone cycle in $\Gamma$, then there is one which is
embedded. To see this, note first that there is one which is immersed (by
\fullref{lm:mc->asc}). Since $b(\Gamma)=1$, the cycle is a covering of
an embedded cycle. Then since the final edge in a monotone cycle appears
only once, the cycle itself must be embedded. 

Hence, we can determine algorithmically whether a given labeled graph
$\Gamma$ with $b(\Gamma) = 1$ contains a strict monotone cycle.
\end{remark}

\begin{proposition}\label{prop:mc=asc}
Suppose $\Gamma\in \LG(G)$ and $b(\Gamma) = 1$. If $\Gamma$ has a strict
monotone cycle and $\Gamma' \in \LG(G)$ then $\Gamma'$ also has a strict
monotone cycle. 
\end{proposition}

\begin{proof}
By \fullref{th:moves} we only need to consider the case when
$\Gamma$ and $\Gamma'$ are related by a slide, induction or $\As^{\pm
1}$--move.  If $\Gamma$ and $\Gamma'$ are related by an induction
move, then both contain strict ascending loops and hence both contain
strict monotone cycles.  Also, if $\Gamma$ and $\Gamma'$ are related
by an $\As^{\pm 1}$--move, then one of the labeled graphs contains a
strict ascending loop and the other contains a strict virtual
ascending loop, hence both contain a strict monotone cycle.  

Now assume that $\Gamma$ has an embedded strict monotone cycle
$(e_{0},\ldots,e_{n},e)$ and that $\Gamma'$ is obtained by sliding an
edge $f$ over an edge $f'$ in $\Gamma$.  Since we can assume that the
strict monotone cycle is embedded, $e_{i} \neq e_{j},\bar{e}_{j}$
for any $i \neq j$.  We have several cases to consider depending on
the configuration of $f,f'$ with respect to the monotone cycle. In all
cases, it suffices to find an edge that can be slid into a loop, since
it will have the same (integral) modulus as $(e_0, \ldots, e_n,
e)$ (because $b(\Gamma) = 1$). 

Clearly if $f,\bar{f} \notin (e_0,\ldots,e_{n},e)$, then this
strict monotone cycle is also a strict monotone cycle in $\Gamma'$.

If $f = e_{i}$ and $f' \neq \bar{e}$, then
$(e_{0},\ldots,e_{i-1},f',e_{i},\ldots,e_{n},e)$ is a strict monotone
cycle in $\Gamma'$.  Likewise, if $f = \bar{e}_{i}$ and $f' \neq e$
then $(e_{0},\ldots,e_{i},\bar{f}',e_{i+1},\ldots,e_{n},e)$ is a
strict monotone cycle in $\Gamma'$.

Since the strict monotone cycle is embedded, the only possible
configurations of $f,f'$ where $f \in \{e_{i},\bar{e}_i\}$ and
$f' \in \{e,\bar{e}\}$ are when $f = e_{0}$ and $f' =
\bar{e}$ or $f = \bar{e}_n$ and $f' = e$.  In the first case
$(e_{1},\ldots,e_{n},e_{0})$ is a strict monotone cycle in $\Gamma'$.
To see this note that $e_{0}$ can slide over $\bar{e}$ and
$\bar{e}$ can slide over $e_{0}$, and hence appropriately chosen
lifts to the Bass--Serre tree carry the same stabilizer.  Then since
$\bar{e}$ can slide over the path $(e_{0},\ldots,e_{n})$, we have that
$\bar{e}_0$ can slide over $(e_{1},\ldots,e_{n})$, after which it becomes
a virtual ascending loop. 
In the second case, with $f = \bar{e}_n$ and $f' = e$, the path
$(e_{0},\ldots,e_{n})$ is a strict monotone cycle in $\Gamma'$ for
similar reasons. 

The remaining cases of interest are when $f \in
\{e,\bar{e}\}$.  If $f = e$, then $(e_{0},\ldots,e_{n},f',e)$ is
a strict monotone cycle in $\Gamma'$, by the following reasoning. For any
strict monotone cycle $(e_0, \ldots, e_n, e)$ we have that 
$\lambda(\bar{e}) q(e_0, \ldots, e_n)$ is an integer, since this is the
label on $\bar{e}$ after sliding over $(e_0, \ldots, e_n)$. Also,
$\lambda(e)$ divides this integer since the modulus of the cycle is
$(\lambda(\bar{e})/\lambda(e)) q(e_0,\ldots, e_n)$. In our situation
$\lambda(f)$ divides $\lambda(e)$, and hence also $\lambda(\bar{e})q(e_0,
\ldots, e_n)$. So in $\Gamma'$ the edge $\bar{e}$ can slide over $(e_0,
\ldots, e_n)$ and then over $f$. 

If $f = \bar{e}$, then $(\bar{f}',e_{0},\ldots,e_{n},e)$ is a strict
monotone cycle in $\Gamma'$, since $\bar{e}$ can slide back over
$\bar{f}'$ and then over $(e_0, \ldots, e_n)$. 
\end{proof}

In the case $b(\Gamma) = 1$, we now have a converse to the first
statement of \fullref{lm:mc->asc}. 

\begin{corollary}\label{co:mc=asc}
If $\Gamma \in \LG(G)$ satisfies $b(\Gamma) = 1$ then $G$ is
ascending if and only if $\Gamma$ has an embedded strict monotone
cycle.
\end{corollary}

\begin{proof}
If $G$ is ascending, then there is a labeled graph $\Gamma' \in 
\LG(G)$ that contains a strict ascending loop (which is a strict monotone
cycle). By \fullref{prop:mc=asc}, $\Gamma$ also contains a strict
monotone cycle. The converse is given by \fullref{lm:mc->asc}. 
\end{proof}

\begin{remark}\label{rm:ascendingcheck}
Note that the latter condition can be checked effectively by 
\fullref{rm:monotonecheck}.  Thus, when $b(\Gamma) = 1$, we can
check algorithmically whether $G$ is ascending.  There is no known
condition for checking whether a GBS group is ascending in
general.
\end{remark}

\begin{definition}\label{def:sGamma}
A mobile edge that is not an ascending loop or the reverse of an
ascending loop is called an \emph{$s$--mobile edge} ($s$ stands for
``slide'').  Note that if $b(\Gamma) = 1$ and $e$ is an ascending loop,
then sliding $e$ or $\bar{e}$ over another edge always results in a graph
that is not reduced. Hence $s$--mobile edges are the only mobile edges
that can slide over another edge while staying inside $\LG(G)$. Given
$\Gamma \in \LG(G)$, let $s(\Gamma)$ be number of geometric $s$--mobile
edges. By the following lemma, this number may also be denoted $s(G)$. 
\end{definition}

\begin{lemma}\label{lm:sGamma}
If $\Gamma, \Gamma' \in \LG(G)$ and $b(\Gamma) = 1$ then $s(\Gamma)
= s(\Gamma')$. 
\end{lemma}

\begin{proof}
As before, we just need to verify this when $\Gamma$ and $\Gamma'$ are
related by a slide, induction or $\As^{\pm 1}$--move.  For slide
moves, the number of mobile edges is invariant (by 
\fullref{co:mobile}) and so is the number of ascending loops.

For the other moves, note that at least one of $\Gamma$ or $\Gamma'$
must be a single strict ascending loop with trees attached.  The
$s$--mobile edges are exactly those which can be slid to and around
the loop.  The result can be verified easily from this description.
\end{proof}

\subsection{The non-mobile subgraph}\label{ssc:nonmobilesubgraph}

Since we are assuming that $b(\Gamma) = 1$ for any $\Gamma \in
\LG(G)$, the image $q(G) \subset \Q^{\times}$ is generated by a single
rational number $q = q(\gamma)$, where $\gamma$ is an (oriented) embedded
cycle in $\Gamma$.  We may assume that $\abs{q} \geq 1$.

Let $\Gamma_{\!  non} \subset \Gamma$ be the \emph{non-mobile
subgraph}, obtained from $\Gamma$ by discarding the mobile edges and any
vertices incident to a strict ascending loop.  Note that
$\Gamma_{\!  non}$ may be disconnected, and may have isolated
vertices.  Let $\Gamma_{1},\ldots,\Gamma_{k}$ be the simply connected
components of $\Gamma_{\!  non}$.  There is at most one component of
$\Gamma_{\!  non}$ not in this list, and this only happens when $G$ is
non-ascending. 

In both of the examples of \fullref{fig:invariant}, the non-mobile
subgraph consists of the two vertices of valence one; all edges are
mobile, and the middle vertex is deleted. 

Each $\Gamma_{i}$ carries a subgroup $G_{i}$ of $G$, well defined up
to conjugacy.  These subgroups and their conjugates will be called
\emph{non-mobile subgroups}.  

For each geometric mobile edge $\{e, \bar{e}\}$, one of its orientations
will be designated as \emph{preferred}. If $e \in \gamma$ then $e$ is
preferred. Otherwise, if $e, \bar{e} \not\in \gamma$, we say $e$ is
preferred if $e$ separates $\iv(e)$ from $\gamma$. Since $b(\Gamma)=1$,
exactly one of $e, \bar{e}$ will have this property. 

\begin{lemma}\label{lm:fi}
For each simply connected component $\Gamma_{i} \subseteq \Gamma_{\!
non}$, there is a unique preferred mobile edge $e_{i}$ such that
$\iv(e_{i}) \in \Gamma_{i}$.  Furthermore, every $s$--mobile edge is
$e_{i}$ or $\bar{e}_i$ for some $i$.  Hence $\Gamma_{\!  non}$ has
exactly $s(G)$ simply connected components. 
\end{lemma}

\begin{proof}
There are two cases depending on whether $\Gamma_{i}$
intersects $\gamma$.  If $\Gamma_{i}$ intersects $\gamma$ (which
can happen if $\gamma$ contains an $s$--mobile edge), then since
$\Gamma_{i}$ is simply connected, there is a (preferred) mobile edge $e_i 
\in \gamma$ such that $\iv(e_i) \in \Gamma_{i}$.  There is at most one
other geometric mobile edge in $\gamma$ that is incident to
$\Gamma_{i}$. Its preferred orientation must meet $\Gamma_{i}$ in its
terminal vertex, since the edges in $\gamma$ are oriented coherently.
For any other preferred mobile edge $e'$ incident to $\Gamma_i$, we have
that $\Gamma_i$ separates $e'$ from $\gamma$, and so $\Gamma_i$ contains
$\tv(e')$, not $\iv(e')$. 

If $\Gamma_{i}$ does not intersect $\gamma$, then since $b(\Gamma) = 1$,
there is a unique preferred mobile edge incident to $\Gamma_{i}$,
separating $\Gamma_i$ from $\gamma$.   

For the second statement, let $e$ be a preferred $s$--mobile edge.  Then
$\iv(e)$ is not the vertex of an ascending loop, and so $\iv(e) \in
\Gamma_{\! non}$. If $\iv(e) \not\in \gamma$ then $\iv(e) \in \Gamma_i$
for some $i$. If $\iv(e) \in \gamma$ then $e\in \gamma$, and hence all
components of $\Gamma_{\! non}$ are simply connected. Thus $\iv(e) \in
\Gamma_i$ for some $i$. In either case, since $\iv(e) \in \Gamma_i$, it
follows that $e = e_i$ by uniqueness. 
\end{proof}

We will be looking carefully at the subgraphs $\Gamma_i$ and how they sit
inside $\Gamma$. For this we need the following definitions. 

\begin{definition}
A \emph{pointed labeled graph} is a triple $\hat\Gamma = (\Gamma, \hat v,
\hat \lambda)$ where $\Gamma$ is a labeled graph, $\hat v \in V(\Gamma)$,
and $\hat \lambda$ is a non-zero integer. It is \emph{reduced} if
$\Gamma$ is reduced and $\hat \lambda \not= \pm 1$. Let $\pLG(G)$
be the set of reduced pointed labeled graphs $(\Gamma, \hat v, \hat
\lambda)$ such that $\Gamma \in \LG(G)$. 

We define an equivalence relation on $\pLG(G)$ via the following
procedure. Given $(\Gamma, \hat v, \hat \lambda)$, adjoin a new
distinguished edge $e$ to $\Gamma$ with $\iv(e) = \hat v$ and label
$\lambda(e) = \hat \lambda$ (the label $\lambda(\bar{e})$ is
irrelevant). Perform any elementary deformation of this graph in which
$e$ is never collapsed. In particular, no edge slides over $e$. Now let
$\hat v' = \iv(e)$, $\hat\lambda' = \lambda(e)$, and delete $e$ to obtain
the labeled graph $\Gamma'$. If $\Gamma'$ is reduced and $\hat \lambda'
\not= \pm 1$, we declare the pointed labeled graphs $(\Gamma, \hat v,
\hat \lambda)$ and $(\Gamma', \hat v', \hat \lambda')$ to be equivalent. 

As always, labeled graphs are considered modulo admissible sign changes,
and this applies to pointed labeled graphs as well. In particular, an
admissible sign change may be performed on the distinguished edge, and so
$(\Gamma, \hat v, \hat \lambda)$ is always equivalent to $(\Gamma, \hat v,
-\hat \lambda)$. Alternatively, this equivalence can be seen by performing
an admissible sign change on every vertex and edge of $\Gamma$. 

Given $\hat\Gamma \in \pLG(G)$ let $\pLG(\hat\Gamma) \subset \pLG(G)$ be
the equivalence class containing $\hat \Gamma$. It is called the
\emph{pointed labeled graph space of $\hat \Gamma$}. 
\end{definition}

\begin{remark}
It is interesting to compare $\pLG(\hat \Gamma)$ with $\LG(\Gamma)$. For
example, let $\hat \Gamma = (\Gamma, \hat v, \hat \lambda)$ where
$\Gamma$ is the labeled graph of \fullref{fig:nomobiles}, $\hat v$ is
the upper left vertex, and $\hat \lambda = 8$. Then $|\pLG(\hat \Gamma)|
= \infty$ even though $|\LG(\Gamma)| < \infty$, because the distinguished
edge can slide around the cycle in the counterclockwise direction,
increasing $\hat \lambda$. Indeed, for any labeled graph $\Gamma$ with a
non-trivial integral modulus, there will be pointed labeled graphs $\hat
\Gamma$ with $| \pLG(\hat \Gamma)| = \infty$, by similar
reasoning. However, if $\Gamma$ has no non-trivial integral moduli, then
we have the following result. 
\end{remark}

\begin{proposition}\label{prop:pointed}
Suppose $G$ has no non-trivial integral moduli. Then
\begin{enumerate}
\item $|\pLG(\hat \Gamma)| < \infty$ for every $\hat \Gamma \in \pLG(G)$,
  and \label{al:finite} 
\item there is an algorithm which, given $\hat \Gamma, \hat \Gamma' \in
  \pLG(G)$, determines whether they are in the same pointed labeled graph
  space. \label{al:algorithm} 
\end{enumerate}
\end{proposition}

\begin{proof} 
Given $\hat \Gamma = (\Gamma, \hat v, \hat \lambda)$ let $\Gamma_0$ be the
reduced labeled graph obtained from $\Gamma$ by adjoining a new edge $e$
and a new vertex $\tv(e)$, with $\iv(e) = \hat v$, $\lambda(e) = \hat
\lambda$, and $\lambda(\bar e) = 2$. Let $G_0$ be the new GBS group. This
operation does not change the image of the modular homomorphism, so $G_0$
has no non-trivial integral moduli. Now observe that $\pLG(\hat \Gamma)$
embeds into $\LG(G_0)$, by identifying the distinguished edge with
$e$. Conclusion \eqref{al:finite} follows because $\LG(G_0)$ is
finite, by \cite[Theorem 8.2]{ar:F06}. 

For \eqref{al:algorithm} one considers elementary deformations of
$\Gamma_0$ in which $e$ is never collapsed. By \fullref{co:slides},
if $\Gamma_0'$ is related to $\Gamma_0$ by such a deformation, then there
is a sequence of slide moves from $\Gamma_0$ to $\Gamma_0'$ in which no 
edge ever slides over $e$. Now, given $\hat \Gamma$ and $\hat \Gamma'$,
start with $\Gamma_0$ and perform all possible sequences of slide moves,
never sliding an edge over $e$. All labeled graphs thus obtained yield
pointed labeled graphs in $\pLG(\hat \Gamma)$ (by recording
$\iv(e)$ and $\lambda(e)$ and deleting $e$). Moreover every pointed
labeled graph in $\pLG(\hat \Gamma)$ will be found, since these slides
take place in $\LG(G_0)$, which is finite. 
\end{proof}

\begin{definition} \label{def:pointed}
Recall that given $\Gamma\in \LG(G)$ with $b(\Gamma) = 1$, each simply
connected component $\Gamma_i$ of $\Gamma_{\! non}$ has a preferred
mobile edge $e_i$ associated to it, with $\iv(e_i) \in \Gamma_i$. Define
$\hat \Gamma_i$ to be the pointed labeled graph $(\Gamma_i, \iv(e_i),
\lambda(e_i))$. This data will also be denoted $(\Gamma_i, \hat v_i, \hat
\lambda_i)$. Note that $\hat\Gamma_i$ is reduced, because $\Gamma$ is,
and so $\hat\Gamma_i \in \pLG(G_i)$ for each $i$. 
\end{definition}

\begin{theorem}\label{th:invariant} 
Suppose $\Gamma, \Gamma' \in \LG(G)$ and $b(\Gamma) = 1$. Then 
\begin{enumerate}
\item $\Gamma$ and $\Gamma'$ define the same non-mobile subgroups of $G$,
  and \label{al:Gi} 
\item for each non-mobile subgroup $G_i$, the corresponding pointed
  labeled graphs $\hat \Gamma_i$ and $\hat \Gamma_i'$ are equivalent in
  $\widehat \LG(G_i)$. \label{al:Gammai} 
\end{enumerate}
\end{theorem}

\begin{proof}
We may assume that $\Gamma$ and $\Gamma'$ are related by a slide,
induction, or $\As$--move. 

First consider an induction move. Both labeled graphs have ascending
loops, where the move takes place, and note that every edge incident to
an ascending loop is mobile. Hence the non-mobile subgraphs and subgroups
do not change, nor do the labels $\lambda(e_i)$ (since $\iv(e_i)$ is not
the vertex of the ascending loop). 

Next suppose that $\Gamma'$ is obtained from $\Gamma$ by an
$\As$--move, exactly as pictured in \fullref{def:A}. The virtually
ascending loop in $\Gamma$ with labels $(k, k\ell m)$ is $e_i$ for some
$i$. Then the vertex of the loop is $\hat v_i$ and $\hat \lambda_i =
k$. After the $\As$--move, the newly created edge with labels $(k, \ell)$
becomes $e_i$, and the subgraph $\Gamma_i$ is unchanged. It is still the
case that $\hat \lambda_i = k$, and $\hat v_i$ has not moved. All other
subgraphs $\Gamma_j$ are also unchanged. Hence $G_i = G_i'$ and $\hat
\Gamma_i = \hat \Gamma_i'$ for all $i$. 

Now suppose that $\Gamma'$ is obtained from $\Gamma$ by sliding $e$ over
$e'$. To prove \eqref{al:Gi} it suffices to show that the simply
connected components of $\Gamma_{\! non}$ contain the same edges and
vertices before and after the slide move. If $e$ is mobile then
$\Gamma_{\! non}$ does not change at all, and \eqref{al:Gi} holds. So
assume that $e$ is non-mobile, which implies that $e'$ is also
non-mobile, by \fullref{lm:nonmobileovermobile}. Now the slide move
takes place entirely within $\Gamma_{\! non}$, and induces a homotopy
equivalence $\Gamma_{\! non} \simeq \Gamma_{\! non}'$ of underlying
graphs. Thus the simply connected components are preserved and
\eqref{al:Gi} holds. 

Now consider part \eqref{al:Gammai}. If $e$ and $e'$ are non-mobile then
the preferred mobile edges $e_i$ do not change, nor do $\iv(e_i)$ and 
$\lambda(e_i)$. Thus $\hat \Gamma_i$ and $\hat \Gamma_i'$ are equivalent
in $\pLG(G_i)$. If $e$ is mobile then it must be an $s$--mobile edge
(cf. \fullref{def:sGamma}) and so $e$ is $e_i$ or $\bar{e}_i$ for
some $i$. If $e=\bar{e}_i$ then $\iv(e_i)$ and $\lambda(e_i)$ do not
change, and $\hat \Gamma_i = \hat \Gamma_i'$ for all $i$. Now suppose
that $e = e_i$. If $e'$ is non-mobile and is in $\Gamma_i$ for some $i$
then $\hat \Gamma_i$ and $\hat \Gamma_i'$ are equivalent in
$\pLG(G_i)$. If $e'$ is non-mobile and not in any $\Gamma_i$ then $\hat
\Gamma_i = \hat \Gamma_i'$ for each $i$. 

Lastly, suppose that $e = e_i$ and $e'$ is mobile. Note that $\iv(e)$ is
not the vertex of an ascending loop, since $e$ is $s$--mobile and
preferred. Hence $e'$ is also an $s$--mobile edge. It is not preferred
because no two preferred mobile edges have a common initial vertex. Thus
$e' = \bar{e}_j$ for some $j \not= i$. Note that before the slide,
$\iv(e) = \iv(e') \in \Gamma_i$ and $\tv(e') \in \Gamma_j$, and after the
slide, $\iv(e) = \tv(e') \in \Gamma_j$ and $\iv(e') \in \Gamma_i$. Thus,
by the uniqueness property of \fullref{lm:fi}, $e_i$ becomes $e_j$ and
$e_j$ becomes $\bar{e}_i$. We also have $\hat v_i = \hat v_i'$ and $\hat
\Gamma_i = \hat \Gamma_i'$ for all $i$. 

The only remaining issue is the labels $\hat \lambda_i, \hat \lambda_j$
and $\hat \lambda_i', \hat \lambda_j'$. We will show that $\lambda(e)
= \pm \lambda(e')$, which implies that $\hat \lambda_i = \pm \hat
\lambda_i'$ and $\hat \lambda_j = \pm \hat \lambda_j'$, completing the
proof. There are two cases. 

If $e', \bar{e}' \not \in \gamma$, then since this is a mobile edge, the
geometric edge $\{e', \bar{e}'\}$ can slide to and around $\gamma$ in the
positive direction. Since $e'$ separates $\tv(e')$ from $\gamma$, the
endpoint $\tv(e')$ can never meet $\gamma$ after sliding $\{e',
\bar{e}'\}$. Hence it is $e'$, and not $\bar{e}'$, which slides to and
around $\gamma$. Such a slide sequence includes a slide of $e'$ over $e$
(whether $e \in \gamma$ or $e \not\in \gamma$). Also, just before this
particular slide, the label $\lambda(e')$ has not changed, since
$\iv(e')$ has remained within a subtree of $\Gamma$ until this
point. Hence $\lambda(e)$ divides $\lambda(e')$. On the other hand, since
$e$ can slide over $e'$, we have that $\lambda(e')$ divides $\lambda(e)$. 

If $\bar{e}' = e_j \in \gamma$ then write $\gamma$ as $(e_i, \gamma_0,
e_j)$. Note that $|\Sl(\Gamma, e')| < \infty$ since $\Gamma - \{e',
\bar{e}'\}$ is a tree. Hence $e'$ is part of a
strict monotone cycle, which we may take to be embedded, and must then be
either $(e_i, \gamma_0, e_j)$ or $(\bar{\gamma}_0, \bar{e}_i,
\bar{e}_j)$. The second case does not occur since this cycle has modulus 
$1/q$, which is not in $\Z - \{\pm 1\}$. So $\bar{e}_j = e'$ can slide
over $e_i = e$, and $\lambda(e)$ divides $\lambda(e')$. But $e$ can slide
over $e'$, and so $\lambda(e')$ divides $\lambda(e)$. 
\end{proof}

\begin{remark}
It can be shown that conclusion \eqref{al:Gi} holds even without the 
assumption that $b(\Gamma) = 1$. More specifically, all three types of
moves preserve the connected components of the non-mobile
subgraph. (Recall from \fullref{rm:nonmobilecanonical} that the set of
non-mobile edges is preserved by the three moves.) 
\end{remark}

\begin{definition} \label{def:plgi} 
We may now define an invariant for non-elementary GBS groups $G$ with
$b(G) = 1$. Choose $\Gamma \in \LG(G)$ and let $\PLGI(G)$ be the collection
of pointed labeled graph spaces $\{\pLG(\hat \Gamma_i) \}$ indexed by the
conjugacy classes of non-mobile subgroups of $G$. By 
\fullref{th:invariant}, $\PLGI(G)$ is independent of the choice of
$\Gamma$. 

Moreover $\PLGI(G)$ is computable: given labeled graphs representing $G$
and $G'$, one may write down representatives for the collections
$\PLGI(G)$ and $\PLGI(G')$, and determine algorithmically whether
$\PLGI(G) = \PLGI(G')$, by \fullref{prop:pointed}. 
\end{definition}

\subsection{Ascending Betti number one GBS groups}\label{ssc:ascending}

Let $G$ be an ascending GBS group with $b(G) = 1$.  Recall
that $q(G) \subset \Q^{\times}$ is generated by $q = q(\gamma)$ where
$\gamma$ is an (oriented) embedded cycle.  Since $G$ is ascending, $q \in
\Z$ and $|q| > 1$.  Let $F(q) \subset \Q^{\times}$ be the subgroup
generated by the integral factors of $q$. We will define an invariant
$\xi(G) \in (\Q^{\times}/ \langle q \rangle)^s / F(q)$, where $s =
s(G)$ and $F(q)$ acts diagonally on the group $(\Q^{\times}/ \langle q
\rangle)^s$. 

Given a labeled graph in $\Gamma \in \LG(G)$, let $\Gamma_i$ and $G_i$ be
defined as in \fullref{ssc:nonmobilesubgraph}, and let $e_1, \ldots,
e_s$ be the preferred mobile edges defined by \fullref{lm:fi}. Also
choose a mobile edge $e \in \gamma$, called the \emph{reference
  edge}. This edge may or may not be among the edges $e_i$, depending on
whether the strict monotone cycle is an ascending loop. Based on $e$, we
will define an element $\xi_i \in \Q^{\times}/\langle q \rangle$ for each
$G_i$, and the resulting $s$--tuple will represent the invariant $\xi(G)$. 

First we claim that there are lifts $\tilde{e}, \tilde{e}_{1},
\ldots,\tilde{e}_{s}$ in the Bass--Serre tree of $\Gamma$ such that
$G_{\tilde{e}_{i}} \subseteq G_{\tilde{e}}$ for each $i$. Note that we
are free to perform slide moves without affecting this claim. If $e$ is a
strict virtually ascending loop, then all mobile edges in $\Gamma$ can be
slid to be adjacent to $e$. Then lifts can be chosen so that
$\tv(\tilde{e}_{i}) = \iv(\tilde{e})$ for each $i$, which implies that
$G_{\tilde{e}_{i}} \subseteq G_{\tilde{e}}$. Otherwise, if $e$ is not a
virtually ascending loop, then it is part of a strict monotone cycle, and
it can be made into a virtually ascending loop by slide moves. Now choose
lifts as before. 

We define $\xi_i = [G_{\tilde{e}}:G_{\tilde{e}_{i}}]$. Note that a
different choice of $\tilde{e}_i$ defines the same element of
$\Q^{\times}/ \langle q \rangle$, because the two lifts are related by an
element of $G$ with modulus a power of $q$. A different choice of
$\tilde{e}$ also makes no difference, by transport of structure. Now
define $\xi(\Gamma) \in (\Q^{\times} / \langle q \rangle )^s / F(q)$ to
be the element represented by $(\xi_1, \ldots, \xi_s)$. 

\begin{lemma}\label{lm:edgeinvariant}
The element $\xi(\Gamma) \in (\Q^{\times}/\langle q \rangle )^s/ F(q)$ is
independent of the choice of reference edge. 
\end{lemma}

\begin{proof}
Consider $\xi'(\Gamma)$ defined using a reference edge $e' \in \gamma$
instead of $e$. We will show that there are lifts $\tilde{e},\tilde{e}'$
such that $G_{\tilde{e}'}\subseteq G_{\tilde{e}}$ and 
$[G_{\tilde{e}} : G_{\tilde{e}'}]$ is a factor of $q$.  Then the
$s$--tuples $(\xi_1, \ldots, \xi_s), (\xi_1', \ldots, \xi_s') \in
(\Q^{\times}/\langle q \rangle)^s$ differ by this factor, and are
equivalent in $(\Q^{\times}/\langle q \rangle)^s/F(q)$. 

Reversing orientations of $e, e'$ if necessary, the cycle $\gamma$ can be
written as $(\alpha,e,\beta,e')$. Both $(\alpha,e,\beta,e')$ and 
$(\beta,e',\alpha,e)$ are strict monotone cycles, because $e$ and $e'$ are
mobile.  Now $\lambda(e)$ divides $q(\alpha)\lambda(\bar{e}')$, as
$(\alpha,e,\beta)$ is an $\bar{e}'$--edge path. Similarly $\lambda(e')$
divides $q(\beta)\lambda(\bar{e})$. Hence the modulus $q$ can be written
as the product of two integers:
\begin{equation*}
q \ = \ \frac{q(\alpha)\lambda(\bar{e}')}{\lambda(e)}
\frac{q(\beta)\lambda(\bar{e})}{\lambda(e')}.
\end{equation*}
Lifting the path $(e', \alpha, e)$ to $(\tilde{e}', \tilde{\alpha},
\tilde{e})$ we obtain $\tilde{e}$ and $\tilde{e}'$ with $G_{\tilde{e}'}
\subseteq G_{\tilde{e}}$. Since $\overline{e}'$ can slide over $(\alpha,
e)$, we have $[G_{\tilde{e}} : G_{\tilde{e}'}] =
\abs{\frac{q(\alpha)\lambda(\bar{e}')}{\lambda(e)}}$. Hence this index
divides $q$. 
\end{proof}

Next we show that $\xi(\Gamma)$ is an invariant of $G$, and hence may be
denoted $\xi(G)$. 

\begin{proposition}\label{th:xi-invariant}
For any two graphs $\Gamma,\Gamma' \in \LG(G)$ we have
$\xi(\Gamma) = \xi(\Gamma')$.
\end{proposition}

\begin{proof}
By \fullref{th:moves}, we may assume that $\Gamma'$ is obtained from
$\Gamma$ by a slide, induction, or $\As^{\pm 1}$--move. We consider the
case of a slide move first. 

Since $b(\Gamma) = 1$, the slide move does not create or remove strict
ascending loops, and so the set of $s$--mobile edges is unchanged. We may
also choose a reference edge $e \in \Gamma$ that remains on the embedded
circuit in $\Gamma'$. Thus, the collection of edges $e, e_1, \ldots, e_s$
and their lifts, used to define $\xi$, can be chosen to agree for
$\Gamma$ and $\Gamma'$. The only change to be accounted for in passing
from $\Gamma$ to $\Gamma'$ is that the correspondence between $s$--mobile
edges and conjugacy classes of non-mobile subgroups may change. That is,
the indexing of the entries of $\xi(\Gamma)$ may change. 

Recall from the proof of \fullref{th:invariant} that if one
$s$--mobile edge slides over another, then their indices and preferred
orientations may be exchanged. However, it was shown that
whenever this occurs, the labels of the two edges at their common vertex
are the same, up to sign. Thus, choosing adjacent lifts $\tilde{e}_i$ and
$\tilde{e}_j$, we have $G_{\tilde{e}_i} = G_{\tilde{e}_j}$, and therefore
$\xi_i = \xi_j$. It follows that $\xi(\Gamma) = \xi(\Gamma')$. 

If $\Gamma$ and $\Gamma'$ differ by an induction move, then there are
strict ascending loops $e \in \Gamma$ and $e' \in \Gamma'$ along which
the move occurs. These edges will be the reference edges for
$\xi$. The $s$--mobile edges for $\Gamma$ and $\Gamma'$ will be the same,
with the same indexing, since the move does not affect the non-mobile
subgraph. Thus we may choose the same lifts $\tilde{e}_i$ for $\Gamma$
and for $\Gamma'$. We may also choose the lifts $\tilde{e}$ and
$\tilde{e}'$ so that $G_{\tilde{e}} \subseteq G_{\tilde{e}'}$ and
$[G_{\tilde{e}'} : G_{\tilde{e}}]$ is a factor $m$ of $q$. (Even though
$\tilde{e}$ and $\tilde{e}'$ are in different trees, this can be
arranged.) Then 
$\xi_i(\Gamma') = m\xi_i(\Gamma)$ for all $i$, and so $\xi(\Gamma') =
\xi(\Gamma)$. 

Now suppose that $\Gamma'$ is obtained from $\Gamma$ by an
$\As^{-1}$--move, exactly as pictured in \fullref{def:A}. In
$\Gamma$, the edge with labels $\ell$ and $k$ is an $s$--mobile edge, say
$e_1$, with initial vertex $v$ on the right. The loop is the reference
edge $e$. Choose a lift $\tilde{e}_1$ and let $\tilde{v}$ be its initial
vertex. The $\As^{-1}$--move does not affect $\tilde{v}$, and the loop
$e_1' \in \Gamma'$ has a lift $\tilde{e}_1'$ with initial vertex
$\tilde{v}$, with the same stabilizer as $\tilde{e}_1$. Note that $e_1'$
is indeed the $s$--mobile edge in $\Gamma'$ corresponding to $G_1$. The
other non-mobile subgraphs and $s$--mobile edges are unchanged. Thus, the
stabilizers of lifts of $s$--mobile edges may be chosen to agree for
$\Gamma$ and $\Gamma'$. What has changed, however, is the reference
edge. The reference edge for $\Gamma'$ is $e_1'$, whose lift
$\tilde{e}_1'$ has stabilizer $G_{\tilde{e}_1}$. The reference edge for
$\Gamma$ is the loop $e$, which has a lift $\tilde{e}$ adjacent to
$\tilde{e}_1$, with $G_{\tilde{e}_1} \subseteq G_{\tilde{e}}$ and
$[G_{\tilde{e}} : G_{\tilde{e}_1}] = \ell$. Now $\xi_i(\Gamma') = \ell
\xi_i(\Gamma)$ for all $i$, and $\xi(\Gamma') = \xi(\Gamma)$, since
$\ell$ divides $q$. 
\end{proof}

Next we define normal forms for the labeled graphs under
discussion. 

\begin{definition}
Suppose $\Gamma$ is a reduced labeled graph with first Betti
number one, in an ascending deformation space. We say that $\Gamma$ is in
\emph{normal form} if it has a strict ascending loop, every mobile edge
is adjacent to this loop, and every label (except possibly the label $q$
on the loop) is positive. Note that if $\Gamma$ is in normal form, then
the $s$--mobile edges are exactly the edges adjacent to the loop, and
$\xi(G)$ is represented by the $s$--tuple $(\lambda(\bar{e}_1), \ldots,
\lambda(\bar{e}_s))$. 

Every $\Gamma$ with $b(\Gamma) = 1$ in an ascending deformation space
can be put into normal form, as follows. First, there is a strict
monotone cycle, which can be made into a strict virtually ascending loop
by slide moves. If necessary, this can be made into a strict ascending
loop by an $\As$--move. Then all $s$--mobile edges can be slid to be
adjacent to the loop. Lastly, since $b(\Gamma)=1$, the labels (other than
$q$) can be made positive by admissible sign changes. 
\end{definition}

\begin{figure}[ht!]
\labellist
\small \hair 2pt
\pinlabel {\scriptsize $2$} [br] at 48 28
\pinlabel {\scriptsize $1$} [tr] at 48 24
\pinlabel {\scriptsize $2$} [br] at 60 31
\pinlabel {\scriptsize $2$} [tr] at 60 21
\pinlabel {\scriptsize $2$} [br] at 91 44
\pinlabel {\scriptsize $2$} [tr] at 91 7
\pinlabel {\scriptsize $2$} [br] at 199 28
\pinlabel {\scriptsize $1$} [tr] at 199 24
\pinlabel {\scriptsize $2$} [br] at 211 31
\pinlabel {\scriptsize $2$} [tr] at 211 21
\pinlabel {\scriptsize $2$} [br] at 242 44
\pinlabel {\scriptsize $4$} [tr] at 242 7
\endlabellist
\begin{center}
\includegraphics{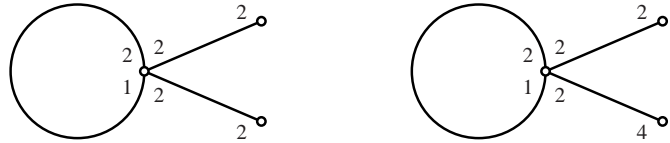}
\end{center}
\caption{Two non-isomorphic GBS groups.} 
\label{fig:invariant}
\end{figure}

\begin{example}
\fullref{fig:invariant} shows two reduced labeled graphs in normal
form representing groups $G$, $G'$. In both cases the invariant $\xi$
is the equivalence class of the pair $(1,1) \in (\Q^{\times}/\langle
2\rangle)^2$. The invariant $\PLGI(G)$ is represented by a pair of
pointed labeled graphs, each consisting of a single vertex, with
distinguished labels $2$ and $2$. On the other hand, $\PLGI(G')$
is represented by two vertices with distinguished labels $2$ and
$4$. Thus, we conclude that $G$ and $G'$ are not isomorphic. Note
that, simple as they are, 
these two groups are not covered by any of the previously known results
on the isomorphism problem (including \fullref{th:iso-finite}). 
\end{example}

\begin{theorem}\label{th:asc-structure} 
Let $G$, $G'$ be ascending Betti number one GBS groups such that
$s(G) = s(G')$ and $q(G) = q(G')$. Then $G$ and $G'$ are isomorphic 
if and only if their non-mobile subgroups are isomorphic, and under this
correspondence between conjugacy classes of non-mobile subgroups, we have 
$\xi(G) = \xi(G')$ and $\PLGI(G) = \PLGI(G')$. 
\end{theorem}

\begin{proof}
Let $(\Gamma,\lambda)$ and $(\Gamma',\lambda')$ be reduced labeled graphs
in normal form representing $G$ and $G'$ respectively. The ``only if''
direction was proved in \fullref{th:invariant} and 
\fullref{th:xi-invariant}. For the other direction we will show that $\Gamma$
and $\Gamma'$ are related by a deformation (considered as unmarked
labeled graphs), which implies that $G \cong G'$. 

Let $G_1, \ldots, G_s$ be the common non-mobile subgroups of $G$ and
$G'$, and let $\Gamma_i$ and $\Gamma_i'$ be the corresponding components
of $\Gamma_{\! non}$ and $\Gamma'_{\! non}$. Then since $\PLGI(G) =
\PLGI(G')$, the pointed labeled graphs $\hat\Gamma_i$ and $\hat\Gamma'_i$
are equivalent in $\pLG(G_i)$. Hence, using the mobile edge $e_i$ as 
the distinguished edge for $\hat\Gamma_i$, there is a deformation of
$\Gamma$, supported in $\hat\Gamma_i \cup e_i$, making $\hat\Gamma_i$
isomorphic to $\hat\Gamma_i'$ as pointed labeled graphs. Thus, we may now
assume that $\hat\Gamma_i$ and $\hat\Gamma_i'$ agree for all $i$. 

Since the graphs are in normal form, the only possible difference
between $\Gamma$ and $\Gamma'$ is in the labels $\lambda(\bar{e}_i)$ and
$\lambda'(\bar{e}_i)$. Since $\xi(G) = \xi(G')$, the $s$--tuples
$(\lambda(\bar{e}_1), \ldots, \lambda(\bar{e}_s))$ and
$(\lambda'(\bar{e}_1), \ldots, \lambda'(\bar{e}_s))$ are equivalent in
$(\Q^{\times}/\langle q \rangle )^s / F(q)$. By performing induction
moves, the $s$--tuples can be made equivalent in $(\Q^{\times}/ \langle q
\rangle )^s$. Now $\lambda(\bar{e}_i)$ and $\lambda'(\bar{e}_i)$ differ
by a factor of a power of $q$. By slide moves of $\bar{e}_i$ over the
ascending loop or its reverse, these labels can be made to agree for
all $i$. 
\end{proof}

\subsection{Non-ascending Betti number one GBS
  groups} \label{ssc:nonascending} 
Let $G$ be a non-ascending GBS group with $b(G) = 1$. Suppose also that
$G$ is not unimodular, and that the modular group $q(G)$ is generated by
an integer $q$. (Otherwise, we are in the situation covered by
\cite{ar:F06}, or alternatively, \fullref{th:iso-finite}.) For now, we
will also assume that $q$ is positive.  
Let the unique embedded cycle $\gamma \subseteq
\Gamma$ be oriented so that $q(\gamma) = q > 1$. 

An edge has infinite slide space if and only if it can 
slide to $\gamma$ and around it at least once in the positive direction.
If it can slide once all the way around, then it can do so infinitely
many times, since its label is multiplied by $q$ each time.  No edge can
slide infinitely many times around in the negative direction, since no
integer is infinitely divisible by $q$. 

Since there are no strict monotone cycles, $\gamma$ does not contain any
mobile edges, and hence is contained in a component $\Gamma_0$ of
$\Gamma_{\! non}$. This is the unique component of $\Gamma_{\!  non}$
that is not simply connected. 

\begin{definition} Let $\Gamma$ be a reduced labeled graph with
$b(\Gamma) = 1$ in a non-ascending, non-unimodular deformation space,
with modulus a positive integer. Let $\gamma \subseteq \Gamma$ be the
unique embedded cycle, oriented so that $q(\gamma) > 1$. We say that
$\Gamma$ is in \emph{normal form} if its labeling is positive and 
every mobile edge is adjacent to $\gamma$, and cannot slide along
$\gamma$ in the negative direction.  Clearly, any $\Gamma$ can be put
into normal form, by sliding the mobile edges to and along $\bar{\gamma}$
as far as they will go. 
\end{definition} 

\begin{theorem}\label{th:normal}
Let $G$ be a non-ascending, non-unimodular GBS group with $b(G) =
1$ and $q(G)$ generated by $q \in \Z_{>0}$. Then $\LG(G)$ contains only
finitely many labeled graphs in normal form, and these can be enumerated
effectively from any $\Gamma \in \LG(G)$. 
\end{theorem}

\begin{proof}
Suppose $\Gamma' \in \LG(G)$ is in normal form.  Let $f_1, \ldots,
f_k$ represent the geometric non-mobile edges of $\Gamma$.  By 
\fullref{co:finitefirst} there are sequences of slide moves \[ \Gamma =
\Gamma^0 \to \Gamma^1 \to \cdots \to \Gamma^k \to \Gamma' \] such that
the moves $\Gamma^{i-1}\to \Gamma^i$ are slides of $f_i, \bar{f}_i$
only, and the moves $\Gamma^k \to \Gamma'$ are slides of mobile edges
only.  Thus we have $\Gamma^i \in \Sl(\Gamma^{i-1},f_i)$ for each $i$,
and since each slide space $\Sl(\Gamma^{i-1},f_i)$ is finite, there
are only finitely many possibilities for the labeled graph $\Gamma^k$.
These graphs can be found effectively by searching the slide spaces
$\Sl(\Gamma^{i-1},f_i)$.  It now suffices to consider the case when
$\Gamma = \Gamma^k$, i.e. when $\Gamma$ and $\Gamma'$ are related by
slide moves of mobile edges only.

The only ambiguity now in determining $\Gamma'$ is in the positioning
and labels of the mobile edges, since the non-mobile subgraphs of
$\Gamma$ and $\Gamma'$ agree.  Note that every mobile edge joins
$\Gamma_0$ to another component $\Gamma_i$ (since $\Gamma$ is in
normal form). Let $G_i$ be the non-mobile subgroup corresponding to
$\Gamma_i$. 

Fix a vertex $v \in \gamma$ and a lift $\tilde{v}$ in the
Bass--Serre tree for $\Gamma$.  Every mobile edge $\bar{e}_i$ may be
slid (in the positive direction) along $\gamma$ to $v$,
after which the label on $\bar{e}_{i}$ is $n_i = [G_{\tilde{v}} :
G_{\tilde{v}} \cap (G_i)^g]$ for some $g \in G$.  Modulo $q$, this
index is independent of $g$, so $[n_i] \in \Q^{\times}/\langle q
\rangle$ depends only on $\Gamma$ and the choice of $v$.

We claim that in fact, $n_i$ itself depends only on the choice of
$v$. Namely, no other representative $q^m n_i$ of $[n_i]$ ($m \in
\Z$) has the property that an edge $e'$ at $v$ with label $q^m n_i$ can
slide around $\gamma$ in the positive direction but not in the negative
direction. To see this, slide the edge with smaller label 
$|m|$ times forward, so the two labels will agree. But then the
other edge could have been slid around $\gamma$ in the negative
direction. 

Now, once $n_i$ is known, the edge $\bar{e}_i$ can be slid back to its
original position in normal form.  This position and the resulting
label on $\bar{e}_i$ are determined by $n_i$. Hence, for any labeled
graph in normal form obtained from $\Gamma$ by sliding mobile edges only,
the labels and initial endpoints of $\bar{e}_{i}$ are uniquely
determined. 

It remains to determine the initial vertices $\iv(e_i)$ and labels
$\lambda'(e_i)$ in $\Gamma'$. The pointed labeled graphs $\hat\Gamma_i$
and $\hat\Gamma_i'$ have the same underlying labeled graphs, and are
equivalent in $\PLGI(G_i)$. Thus, all possible initial vertices
$\hat{v}_i' = \iv(e_i)$ and labels $\hat\lambda_i' = \lambda'(e_i)$ are
obtained by sliding the initial endpoint of $e_i$ within $\Gamma_i$, by
\fullref{co:slides}.  Since $\Gamma_i$ is simply connected, this
slide space is finite and can be searched effectively (cf. 
\fullref{prop:pointed}). 
\end{proof}

We can now prove \fullref{th:iso}. Recall that this theorem solves
the isomorphism problem in the case where one of the labeled graphs has
first Betti number at most one. 

\begin{proof}[Proof of \fullref{th:iso}] 
Let $\Gamma$ and $\Gamma'$ be labeled graphs defining GBS groups $G$ and
$G'$, where $b(\Gamma)\leq 1$. If $q(G)$ is not generated by an integer
then the algorithm of \cite[Corollary 8.3]{ar:F06} determines whether $G
\cong G'$. Hence we may assume that $b(\Gamma) = 1$ and $q(G)$ is
generated by $q \in \Z$ with $\abs{q}>1$. We may also assume that
$b(\Gamma') = 1$ and $q(G') = q(G)$, since otherwise $G \not\cong G'$. 
Moveover, we may assume that $q$ is positive, by 
\fullref{lm:signs}, since the orientation homomorphisms of $\Gamma$ and
$\Gamma'$ agree. 

Now make both graphs reduced by performing collapse moves, and check
whether $\Gamma$ and $\Gamma'$ are ascending (cf. 
\fullref{rm:ascendingcheck}). If one is ascending and the other is not, the
groups are not isomorphic. If both are ascending, then put both into
normal form and verify that $s(\Gamma) = s(\Gamma')$ (if not, then $G
\not\cong G'$). Then identify the subgraphs $\Gamma_i$, $\Gamma_i'$ and
consider permutations $\sigma \in S_s$. For each permutation, check
whether $G_i \cong G_{\sigma(i)}'$ for all $i$ (these GBS groups are
unimodular, so they can be compared). If so, call $\sigma$ an
\emph{admissible permutation} and then re-index the components of
$\Gamma_{\! non}'$ using $\sigma$, so that $G_i
\cong G_i'$ for all $i$. Evaluate and compare the invariants $\xi(G),
\xi(G')$ and $\PLGI(G), \PLGI(G')$, using 
\fullref{prop:pointed}. By \fullref{th:asc-structure}, $G$ 
and $G'$ are isomorphic if these invariants agree. If the invariants
disagree for every admissible permutation, then $G \not\cong G'$, again
by \fullref{th:asc-structure}. 

If both graphs are non-ascending, then put them into normal form. Using
\fullref{th:normal}, enumerate from $\Gamma$ all labeled graphs in
$\LG(G)$ in normal form. Then $G \cong G'$ if and only if $\Gamma'$ is on
this list. 
\end{proof}

%
%
%
%


\def\cprime{$'$}

\end{document}